\PassOptionsToPackage{dvipsnames,table}{xcolor}

\documentclass[10pt,a4paper]{article}
\usepackage[margin=3.24cm]{geometry}
\usepackage{amsfonts,amsmath,bm,dsfont,mathrsfs}
\usepackage{amssymb,amsthm}
\usepackage{enumitem,cases}

\usepackage{tabularx} 
\usepackage{pdflscape}  
\usepackage{rotating}
\usepackage{graphicx}
\usepackage{subcaption}

\usepackage{epstopdf}
\usepackage{algorithm}
\usepackage{algpseudocode}
\usepackage{tikz}
\usepackage{booktabs}    
\usepackage{float}   
\usepackage{multirow}

\usepackage{natbib}
\usepackage[dvipsnames,table]{xcolor}

\usepackage{hyperref}[6.83]
\definecolor{refblue}{RGB}{26,13,171}
\hypersetup{
  colorlinks = true,
  allcolors = refblue,
  urlcolor = BrickRed,
}

\newtheorem{assumption}{Assumption}

\newtheorem{theorem}{Theorem}
\newtheorem{proposition}{Proposition}

\newtheorem{remark}{Remark}
\numberwithin{equation}{section}
\numberwithin{lemma}{section}
\numberwithin{example}{section}
\numberwithin{definition}{section}
\numberwithin{assumption}{section}
\numberwithin{theorem}{section}
\numberwithin{proposition}{section}
\numberwithin{corollary}{section}
\numberwithin{remark}{section}

\newcommand{\email}[1]{\href{mailto:#1}{\texttt{#1}}}

\setlist[description]{style=multiline,leftmargin=2em}
\setlist[itemize]{style=standard,leftmargin=2em}
\setlist[enumerate]{style=standard,leftmargin=2.2em,itemsep=3pt}

\title{HPR-QP: A dual  Halpern Peaceman--Rachford method for solving large-scale convex composite quadratic programming\thanks{The work of Defeng Sun was supported by the Research Center for Intelligent Operations Research, RGC {Senior Research Fellow Scheme No. SRFS2223-5S02}, and  {GRF Project No. 15307822}. The work of Yancheng Yuan was supported by the Research Center for Intelligent Operations Research. The work of Xinyuan Zhao was supported in part by the National Natural Science Foundation of China under Project No. 12271015.
}}
\author{Kaihuang Chen\thanks{Department of Applied Mathematics, The Hong Kong Polytechnic University, Hung Hom, Hong Kong,
  (\email{kaihuang.chen@connect.polyu.hk}).}
\and Defeng Sun\thanks{Department of Applied Mathematics, The Hong Kong Polytechnic University, Hung Hom, Hong Kong,
(\email{defeng.sun@polyu.edu.hk}).} \and 
Yancheng Yuan\thanks{Department of Applied Mathematics, The Hong Kong Polytechnic University, Hung Hom, Hong Kong,
(\email{yancheng.yuan@polyu.edu.hk}).}
\and Guojun Zhang\thanks{Department of Applied Mathematics, The Hong Kong Polytechnic University, Hung Hom, Hong Kong,
(\email{guojun.zhang@connect.polyu.hk}).}
\and Xinyuan Zhao\thanks{Department of Mathematics, Beijing University of Technology, Beijing, P.R. China,
(\email{xyzhao@bjut.edu.cn}).}
}

\def\cK{{\cal K}}

\def\cT{{\cal T}}
\def\cS{{\cal S}}

\def\hx{{\hat x}}

\def\by{{\bar y}}
\def\bz{{\bar z}}
\def\bs{{\bar s}}
\def\bx{{\bar x}}

\def\bw{{\bar w}}
\def\bv{{\bar v}}

\def\bu{{\bar u}}

\begin{document}
\maketitle

\begin{abstract}
In this paper, we introduce {HPR-QP}, a dual Halpern Peaceman--Rachford (HPR) method designed for solving large-scale convex composite quadratic programming. One distinctive feature of HPR-QP is that, instead of working with the primal formulations,  it builds on the novel restricted Wolfe dual introduced in recent years. It also leverages the symmetric Gauss–Seidel technique to simplify subproblem updates without introducing auxiliary slack variables that typically lead to slow convergence. By restricting updates to the range space of the Hessian of the quadratic objective function, HPR-QP employs proximal operators of smaller spectral norms to speed up the convergence. Shadow sequences are elaborately constructed to deal with the range space constraints. Additionally, {HPR-QP} incorporates adaptive restart and penalty parameter update strategies, derived from the HPR method’s \({O}(1/k)\) convergence in terms of the Karush--Kuhn--Tucker residual, to further enhance its performance and robustness. Extensive numerical experiments on benchmark data sets using a GPU demonstrate that our Julia implementation of {HPR-QP} significantly outperforms state-of-the-art solvers in both speed and scalability.


\medskip
\noindent
{\bf MSCcodes:}
90C20, 90C06, 	90C25, 65Y20
\end{abstract}

\section{Introduction}\label{sec:Intro}
In this paper, we develop a dual Halpern Peaceman--Rachford (HPR) method with semi‑proximal terms~\citep{sun2025accelerating} for solving the large-scale convex composite quadratic programming (CCQP) problem:
\begin{equation} \label{QP-primal}
 \min_{x \in \mathbb{R}^n} \left\{ \displaystyle \frac{1}{2} \langle x, Q x \rangle + \langle c, x \rangle + \phi(x) \,\middle|\, A x \in \mathcal{K} \right\},
\end{equation}
where \( Q:\mathbb{R}^{n}\to\mathbb{R}^{n}\) is a self-adjoint positive semidefinite linear operator, \( c \in \mathbb{R}^n \) is a given vector, and \( \phi : \mathbb{R}^n \rightarrow (-\infty,+\infty] \) is a proper, closed, and convex function. Here, \( A: \mathbb{R}^n \to \mathbb{R}^{m} \) is a linear operator, and \(\mathcal{K} := \{ y \in \mathbb{R}^m \mid -\infty \leq l_i \leq y_i \leq u_i \leq +\infty,  1 \leq i \leq m\}\) is a simple polyhedral set. A key feature of our approach is that it does not require an explicit matrix representation of \( Q \), which makes the proposed method particularly suitable for large-scale or matrix-free settings—e.g., when \( Q \) is defined implicitly via Kronecker products or structured operators \citep{anstreicher2001new}. In particular, CCQP includes the classical convex QP (CQP):
\begin{equation} \label{QP-primal-0}
\min_{x \in \mathbb{R}^n} \left\{ \displaystyle\frac{1}{2} \langle x, Q x \rangle + \langle c, x \rangle + \delta_{\mathcal{C}}(x) \,\middle|\, A x \in \mathcal{K} \right\}
\end{equation}
as an important special case, where \(\delta_{\mathcal{C}}(\cdot)\) is the indicator function of the box constraint set \(\mathcal{C} = \{ x \in \mathbb{R}^n \mid L \leq x \leq U \}\), with \(L \in (\mathbb{R} \cup \{-\infty\})^n\) and \(U \in (\mathbb{R} \cup \{+\infty\})^n\).

For solving the CQP \eqref{QP-primal-0}, commercial solvers such as Gurobi~\citep{gurobi} and CPLEX~\citep{cplex2009v12} typically rely on active set methods or interior-point methods. While these methods are robust and effective for small- to medium-sized problems, they face significant scalability challenges when applied to large-scale instances. In particular, interior-point methods suffer from high per-iteration computational costs and substantial memory demands. On the other hand, active set methods—while benefiting from cheaper per-iteration costs—are inherently sequential, making them difficult to parallelize. Even GPU-accelerated interior-point implementations such as CuClarabel~\citep{goulart2024clarabel,chen2024cuclarabel}, which exploit mixed-precision arithmetic, still rely on direct factorization routines, limiting their scalability on large-scale problems.

To overcome these limitations, a variety of first-order solvers, such as SCS~\citep{o2016conic,o2021operator},  OSQP~\citep{stellato2020osqp}, PDQP \citep{lu2023practical}, ABIP \citep{lin2021admm,deng2024enhanced}, PDHCG \citep{huang2024restarted}, and PDCS \citep{lin2025pdcs},  has been developed for large-scale CQP problems.  In particular, SCS~\citep{o2016conic,o2021operator}  applies the Douglas--Rachford method~\citep{lions1979splitting} to solve convex conic programming with convex quadratic objective functions, while OSQP~\citep{stellato2020osqp} implements a generalized ADMM~\citep{eckstein1992douglas} tailored for CQP problems. Both solvers support indirect methods for solving the linear systems that arise at each iteration, which improves scalability over direct factorization routines. Recently, Lu and Yang~\citep{lu2023practical} proposed PDQP, which combines the accelerated primal-dual hybrid gradient method~\citep{chenOptimal2014} with adaptive step sizes and restart strategies—algorithmic enhancements used in the award-winning solver PDLP~\citep{applegate2021practical,applegate2023faster,lu2023cupdlp}. Notably, each step of PDQP has explicit update formulas, eliminating the need for linear system solves. Experiments in~\citep{lu2023practical} show that PDQP, implemented in Julia with GPU support, outperforms SCS (on GPU and CPU) and OSQP (on CPU) on large-scale synthetic CQP problems.

Beyond these developments, some theoretical and algorithmic advances have been made in accelerating the (semi-proximal) Peaceman--Rachford (PR) method \citep{zhang2022efficient,yang2025accelerated,sun2025accelerating,zhang2025hot} using the Halpern iteration~\citep{halpern1967fixed,sabach2017first,lieder2021convergence}. Particularly, Zhang et al. \cite{zhang2022efficient} developed the HPR method without proximal terms by applying the Halpern iteration to the PR method \cite{eckstein1992douglas,lions1979splitting}, achieving an ${O}(1/k)$ iteration complexity in terms of the Karush–Kuhn–Tucker (KKT) residual and the objective error. Sun et al.~\citep{sun2025accelerating} reformulated the semi-proximal PR method as a degenerate proximal point method (dPPM)~\citep{bredies2022degenerate} with a positive semidefinite preconditioner, and applied the Halpern iteration to derive the HPR method with semi-proximal terms, which also enjoys an \({O}(1/k)\) iteration complexity.  Building on this, Chen et al.\citep{chen2024hpr} introduced HPR-LP, a GPU-accelerated solver for large-scale linear programming (LP), which demonstrated significantly better performance than the award-winning solver PDLP \citep{applegate2021practical,applegate2023faster,lu2023cupdlp}. Given PDQP’s strong performance on CQP problems as an extension of PDLP, and HPR-LP’s superior results over PDLP on LP tasks, we are motivated to develop a GPU-accelerated HPR method for solving large-scale CCQP problems, including CQP problems.

A natural approach is to apply the HPR method directly to the primal CCQP \eqref{QP-primal} by introducing a single auxiliary slack variable \(s\):
\begin{equation}\label{QP-primal-s} 
\min_{(x,s)\in \mathbb{R}^{n}\times \mathbb{R}^{m}}\left\{\displaystyle  \frac{1}{2}\langle x,Qx\rangle + \langle c,x\rangle + \phi(x) + \delta_{\mathcal{K}}(s)\;\middle|\;Ax=s\right\}.
\end{equation}
To simplify the subproblem solving within the HPR framework, this approach typically requires a large proximal operator of the form
\[
\displaystyle \mathcal{S}_x = \lambda_Q I_n - Q + \sigma\bigl(\lambda_A I_n - A^*A\bigr),
\]
where \(\lambda_Q \geq \lambda_1(Q)\) and \(\lambda_A \geq \lambda_1(A^*A)\) are constants ensuring that \(\mathcal{S}_x\) is positive semidefinite, and \(\sigma > 0\) is a penalty parameter; \(\lambda_1(\cdot)\) denotes the largest eigenvalue of a self-adjoint linear operator. However, the resulting large spectral norm \(\|\sqrt{\mathcal{S}_x}\|\) can significantly slow convergence (see Appendix~\ref{appendix-primal-1} for algorithmic details). To decouple \(A\) and \(Q\) (so that each can be handled separately), one may introduce a second auxiliary variable \(v\):
\begin{equation}\label{QP-primal-sv}
\min_{(x,s,v)\in \mathbb{R}^{n}\times \mathbb{R}^{m}\times \mathbb{R}^{n}}\left\{\displaystyle
\frac{1}{2}\langle v, Q v \rangle + \langle c, x \rangle + \phi(x) + \delta_{\mathcal{K}}(s)
\;\middle|\; Ax = s,\; x = v
\right\}.
\end{equation}
This yields two simpler proximal operators:
\[
\displaystyle \mathcal{S}_v = \lambda_Q I_n - Q, \qquad \mathcal{S}_x = \sigma(\lambda_A I_n - A^*A),
\]
but the combined spectral radius \(\|\sqrt{\cS_v}\| + \|\sqrt{\cS_x}\|\) remains large, limiting convergence speed (see Appendix~\ref{appendix-primal-2} for algorithmic details).

In this work,  instead of working on the primal forms, we pay our attention to a   novel restricted Wolfe dual of problem~\eqref{QP-primal}, as recently introduced in~\citep{li2018qsdpnal}:

\begin{equation}\label{QP-dual}
\min_{(y,w,z)\in\mathbb{R}^m\times\mathbb{R}^n\times\mathbb{R}^n}
\left\{
\frac{1}{2}\langle w, Qw\rangle + \delta^*_{\mathcal{K}}(-y) + \phi^*(-z)
\;\Big|\; -Qw + A^*y + z = c,\quad w \in \mathcal{W}
\right\},
\end{equation}
where $\mathcal{W}:=\operatorname{Range}(Q)$, the range space of $Q$,  is explicitly imposed in the constraints, as opposed to the classical Wolfe dual form \citep{DornQPDual60,WolfeDual1961} with $\mathcal{W}$ to be taken as the whole space $\mathbb{R}^n$. Since \(Q\) is positive definite on \(\mathcal{W}\), one may consider applying a convergent three-block semi-proximal ADMM~\citep{li2015convergent} to solve problem~\eqref{QP-dual}. However, its convergence guarantee requires choosing the penalty parameter \(\sigma\) proportional to the smallest positive eigenvalue of \(Q\), which leads to slow practical performance. To overcome this limitation, we propose {HPR-QP}, a dual HPR method to solve problem~\eqref{QP-dual}. The main features of {HPR-QP} are summarized below:
\begin{enumerate}
  \item By leveraging the symmetric Gauss--Seidel (sGS) technique~\citep{li2016schur,li2019block}, HPR-QP decouples operators \(A\) and  \(Q\) without introducing auxiliary slack variables. Furthermore, restricting updates to \(\mathcal{W}\) allows {HPR-QP} to employ proximal operators with significantly smaller spectral norms for handling \(Q\), thereby accelerating convergence. Moreover, shadow sequences are constructed to address the numerical challenges introduced by subspace constraints.  Together, these innovations eliminate the need for the large proximal terms required by primal formulations, offering both theoretical and computational advantages. Numerical results in Section~\ref{sec:4} confirm that our restricted Wolfe dual approach substantially outperforms its primal counterparts.  
  \item HPR-QP incorporates adaptive restart and penalty parameter update strategies, derived from the HPR method’s \({O}(1/k)\) iteration complexity in terms of the KKT residual, to further enhance its performance and robustness. Also, {HPR-QP} does not require an explicit matrix representation of \(Q\), allowing it to handle extremely large-scale problem instances. For example, it can efficiently solve  CQP relaxations of quadratic assignment problems (QAPs) involving up to 8,192 locations.
  \item Extensive numerical experiments on benchmark data sets using a GPU demonstrate that our Julia implementation of {HPR-QP} significantly outperforms state-of-the-art solvers in both speed and scalability.
\end{enumerate}

The remainder of this paper is organized as follows. Section~\ref{sec:2} introduces the HPR method with semi-proximal terms for solving the restricted Wolfe dual of CCQP. Section~\ref{sec:3} details the implementation of {HPR-QP}, including its adaptive restart strategy and penalty parameter update scheme. In Section~\ref{sec:4}, we present extensive numerical results across various CCQP benchmark data sets. Finally, Section~\ref{sec:5} concludes the paper with a summary and future directions.

\paragraph{Notation.}

Let \(\mathbb{R}^n\) denote the \(n\)-dimensional real Euclidean space, equipped with the standard inner product \(\langle \cdot, \cdot \rangle\) and its induced norm \(\|\cdot\|\). The infinity norm is denoted by \(\|\cdot\|_\infty\), and the nonnegative orthant is denoted by \(\mathbb{R}_+^n\). For a linear operator  \(A:\mathbb{R}^{n} \to  \mathbb{R}^m\), we denote its adjoint by \(A^*\), and its spectral norm by \(\|A\| := \sqrt{\lambda_1(AA^*)}\), where \(\lambda_1(\cdot)\) denotes the largest eigenvalue of a self-adjoint linear operator. Furthermore, for any self-adjoint positive semidefinite linear operator \(\mathcal{M}: \mathbb{R}^n \to \mathbb{R}^n\), we define the semi-norm \(\|x\|_{\mathcal{M}} := \sqrt{\langle x, \mathcal{M} x \rangle}\) for any \(x \in \mathbb{R}^n\). Given a convex function \(f: \mathbb{R}^n \to (-\infty, +\infty]\), we denote its effective domain by \(\operatorname{dom}(f) := \{x \in \mathbb{R}^n \mid f(x) < +\infty\}\), its subdifferential by \(\partial f(\cdot)\), its convex conjugate by \(f^*(z) := \sup_{x \in \mathbb{R}^n} \{\langle x, z \rangle - f(x)\}\), and its proximal mapping by \(\operatorname{Prox}_f(x) := \arg\min_{z \in \mathbb{R}^n} \left\{ f(z) + \frac{1}{2} \|z - x\|^2 \right\}\), respectively. Moreover, let \(C \subseteq \mathbb{R}^n\) be a closed convex set. The indicator function  over \(C\), denoted by \(\delta_C(\cdot)\),     is defined as  \(\delta_C(x)=0\)  if \(x \in C\) and \(\delta_C(x)=+\infty\) if \(x \notin C\). The Euclidean distance from a point \(x \in \mathbb{R}^n\) to \(C\) is \(\operatorname{dist}(x, C) := \inf_{z \in C} \|z - x\|\), and the Euclidean projection onto \(C\) is \(\Pi_C(x) := \arg\min_{z \in C} \|x - z\|\).

\section{A Dual HPR Method for Solving CCQP}\label{sec:2}

In this section, we present a general HPR framework with semi-proximal terms for solving the restricted Wolfe dual problem~\eqref{QP-dual}. In particular, it includes a variant based on the sGS technique for simplifying the solution of the subproblems as a special instance.

\subsection{An HPR Method with Semi-proximal Terms}
According to \cite[Corollary 28.3.1]{rockafellar1970convex}, a point \( (y^*, w^*, z^*) \in \mathbb{R}^{m} \times \mathcal{W} \times \mathbb{R}^{n} \) is an optimal solution to problem~\eqref{QP-dual} if and only if there exists \( x^* \in \mathbb{R}^n \) such that the KKT system below is satisfied:
\begin{equation}\label{eq:KKT-QP}
	Q w^* - Q x^* = 0, \quad A x^* \in \partial \delta^*_{\mathcal{K}}(-y^*), \quad x^* \in \partial \phi^*(-z^*), \quad -Qw^{*}+A^* y^* + z^* - c = 0.
\end{equation}
Let  \( \sigma > 0 \) be a given penalty parameter.  Define the augmented Lagrangian function \( L_\sigma(y, w, z; x) \) associated with problem~\eqref{QP-dual}, for any \( (y, w, z, x) \in \mathbb{R}^m \times \mathcal{W} \times \mathbb{R}^n \times \mathbb{R}^n \), as follows
$$
L_\sigma(y,w,z ; x)=\frac{1}{2}\langle w, Q w\rangle+ \delta_{\cK}^{*}(-y)+\phi^{*}(-z) + \langle x,  -Qw+A^{*}y+z-c \rangle  +\frac{\sigma}{2}\left\|-Q w+A^*y+z-c\right\|^2.
$$
For notational convenience, we denote the tuple \( (y, w, z, x) \) by \( u \), and define the space \( \mathcal{U} := \mathbb{R}^{m} \times \mathcal{W} \times \mathbb{R}^{n} \times \mathbb{R}^n \). The HPR method with semi-proximal terms, which corresponds to the accelerated preconditioned ADMM with \( \alpha = 2 \) proposed in \citep{sun2025accelerating}, is presented in Algorithm~\ref{alg:acc-pADMM-QP-0} for solving the restricted Wolfe dual problem~\eqref{QP-dual}.

\begin{algorithm}[H] 
	\caption{An HPR method for solving the restricted Wolfe dual problem~\eqref{QP-dual}}
	\label{alg:acc-pADMM-QP-0}
	\begin{algorithmic}[1] 
		\State \textbf{Input:} Choose a self-adjoint positive semidefinite linear operator \( \mathcal{T}_1 \) on \( \mathbb{R}^{m} \times \mathcal{W} \). Denote $u=(y,w,z,x)$ and $\bu=(\by,\bw,\bz,\bx)$. Let \( u^0 = (y^0, w^0, z^0, x^0) \in \mathcal{U} \), and set  \( \sigma > 0 \).
		\For{\( k = 0, 1, 2, \ldots \)}
		      \State {Step 1. $\displaystyle \bz^{k+1}=\underset{z \in \mathbb{R}^{n}}{\arg \min }\left\{L_\sigma(y^k, w^k,z; x^k)\right\} $;} 
                \State {Step 2. $\displaystyle  \bx^{k+1}={x}^k+\sigma (-Qw^{k}+A^{*}y^{k}+\bz^{k+1}-c)$;}
                \State {Step 3. $\displaystyle (\by^{k+1},\bw^{k+1})=\underset{(y,w) \in \mathbb{R}^{m}\times \mathcal{W}}{\arg \min }\left\{L_\sigma(y,w,\bz^{k+1}; \bx^{k+1})+\frac{1}{2}\|(y,w)-(y^k,w^k)\|^2_{\cT_1}\right\}$;}      
                \State {Step 4. $\displaystyle \hat{u}^{k+1}= 2 \bu^{k+1}-u^k$;}
	          \State {Step 5. $\displaystyle u^{k+1}=\frac{1}{k+2}u^{0}+\frac{k+1}{k+2}\hat{u}^{k+1}$;}
		\EndFor
	\end{algorithmic}
\end{algorithm}

To analyze the global convergence of the HPR method with semi-proximal terms presented in Algorithm~\ref{alg:acc-pADMM-QP-0}, we define the linear operator \( \mathcal{H}: \mathbb{R}^{m} \times \mathcal{W} \to \mathbb{R}^{m} \times \mathcal{W} \) as follows:
\begin{equation}\label{Hessian}
\mathcal{H} := \sigma A_Q^{*}A_Q + 
\begin{pmatrix}
\mathbf{0} & \mathbf{0} \\
\mathbf{0} & Q
\end{pmatrix}
=
\begin{pmatrix}
\sigma A A^{*} & -\sigma A Q \\
- \sigma Q A^{*} & \sigma Q^2 + Q
\end{pmatrix},
\end{equation}
where  \( A_Q:=[A^{*} -Q]\in \mathbb{R}^{n\times (m+n)}\). Furthermore, we define a self-adjoint linear operator \( \mathcal{M} : \mathbb{R}^{m} \times \mathcal{W} \times \mathbb{R}^{n} \times \mathbb{R}^{n} \to \mathbb{R}^{m} \times \mathcal{W} \times \mathbb{R}^{n} \times \mathbb{R}^{n} \) as
\begin{equation}\label{def:M-QP}
\mathcal{M} =\begin{bmatrix}
\sigma A_Q^{*}A_Q
+ \mathcal{T}_1 & \quad  \mathbf{0}& \quad 
A_{Q}^{*} \\
\mathbf{0}& \quad \mathbf{0} & \quad  \mathbf{0}&\\
A_{Q} & \quad \mathbf{0} &\quad  \frac{1}{\sigma} I_n
\end{bmatrix},
\end{equation}
where \( I_n \in \mathbb{R}^{n \times n} \) denotes the identity matrix. Now, we make the following assumptions:
\begin{assumption}\label{ass: CQ-QP}
There exists a vector $({y}^*,w^{*}, {z}^*,{x}^*)\in \mathbb{R}^{m}\times \mathcal{W}\times \mathbb{R}^{n}\times \mathbb{R}^{n}$ satisfying the KKT system \eqref{eq:KKT-QP}.
\end{assumption}
\begin{assumption}\label{ass: stronglyconvex}
The operator $\cT_1$ is a self-adjoint positive semidefinite linear operator such that $ \cT_1+\mathcal{H}$ is positive definite on $\mathbb{R}^{m}\times \mathcal{W}$.      
\end{assumption}
Under Assumption \ref{ass: CQ-QP}, solving problems
\eqref{QP-primal} and \eqref{QP-dual} is equivalent to finding a $u^{*} \in \mathcal{U}$ such that $\mathbf{0}\in \cT u^{*}$, where the maximal monotone operator $\cT$ is defined by
	\begin{equation}\label{def:T-QP}
		\cT u=\left(\begin{array}{c}
-\partial \delta_{\mathcal{K}}^{*}(-y)+Ax\\
Qw -Qx\\
-\partial \phi^{*}(-z)+x\\
c-A^{*}y+Qw-z
		\end{array}     \right) \quad  \forall u=(y,w,z,x)\in \mathbb{R}^{m} \times  \mathcal{W} \times \mathbb{R}^{n} \times \mathbb{R}^{n}.
	\end{equation}
Moreover, under Assumption~\ref{ass: stronglyconvex}, each subproblem in Algorithm~\ref{alg:acc-pADMM-QP-0} admits a unique solution. Based on Corollary 3.5 in \citep{sun2025accelerating}, we establish the following global convergence result for the HPR method.
\begin{proposition}\label{Prop:convergence}
Suppose that Assumptions \ref{ass: CQ-QP} and \ref{ass: stronglyconvex} hold. Then the sequence \(\{\bu^k\} = \{(\by^k, \bw^k,\bz^k, \bx^k)\}\) generated by the HPR method with semi-proximal terms in Algorithm \ref{alg:acc-pADMM-QP-0} converges to the point \(u^{*} = (y^*, w^{*}, z^*, x^*)\), where \((y^*,w^{*}, z^*)\) solves problem \eqref{QP-dual} and \(x^*\) solves problem \eqref{QP-primal}.
\end{proposition}

To further analyze the complexity of the HPR method with semi-proximal terms in terms of the KKT residual and the objective error, we consider the residual mapping associated with the KKT system~\eqref{eq:KKT-QP}, as introduced in~\citep{han2018linear}:
\begin{equation}\label{def:KKT_residual}
\mathcal{R}(u) = 
\begin{pmatrix}
Ax - \Pi_{\mathcal{K}}(Ax  - y) \\
Qw - Qx \\
x - \operatorname{Prox}_{\phi}(x - z) \\
c - A^{*}y + Qw - z
\end{pmatrix} \quad \forall\, u = (y, w, z, x) \in \mathbb{R}^{m} \times \mathcal{W} \times \mathbb{R}^{n} \times \mathbb{R}^{n}.
\end{equation}
In addition, let \( \{(\by^k, \bw^k, \bz^k)\} \) be the sequence generated by Algorithm~\ref{alg:acc-pADMM-QP-0}. We define the objective error as
\[
\begin{aligned}
h(\by^{k+1}, \bw^{k+1}, \bz^{k+1}) := &\ \frac{1}{2} \langle \bw^{k+1}, Q \bw^{k+1} \rangle + \delta_{\mathcal{K}}^*(-\by^{k+1})  + \phi^*(-\bz^{k+1}) \\
& - \left( \frac{1}{2} \langle w^*, Q w^* \rangle + \delta_{\mathcal{K}}^*(-y^*)  + \phi^*(-z^*) \right) \quad \forall\, k \geq 0,
\end{aligned}
\]
where \( (y^*, w^*, z^*) \) is a solution to the dual problem \eqref{QP-dual}. Based on the iteration complexity results in Proposition 2.9, Theorem 3.7, and Remark 3.8 of~\citep{sun2025accelerating}, we obtain the following iteration complexity bounds for the HPR method with semi-proximal terms.

\begin{proposition}\label{prop:complexity-acc-pADMM}
  Suppose that Assumptions~\ref{ass: CQ-QP} and~\ref{ass: stronglyconvex} hold. Let \( \{\bu^k\} = \{(\by^k, \bw^k, \bz^k, \bx^k)\} \) and \( \{u^k\} = \{(y^k, w^k, z^k, x^k)\} \) be the sequences generated by the HPR method with semi-proximal terms in Algorithm~\ref{alg:acc-pADMM-QP-0}, and let \( u^* = (y^*, w^*, z^*, x^*) \) be a solution to the KKT system~\eqref{eq:KKT-QP}. Define \( R_0 := \|u^0 - u^*\|_{\mathcal{M}} \). Then, for all \( k \geq 0 \), the following complexity bounds hold:
\begin{equation*}\label{eq: complexity-M}
\|\bu^{k+1} - u^k\|_{\mathcal{M}} \leq \frac{R_0}{k+1},
\end{equation*}
\begin{equation*}\label{eq: complexity-bound-KKT}
\|\mathcal{R}(\bar{u}^{k+1})\| \leq \left( \frac{\sigma \|A_{Q}^{*}\| + 1}{\sqrt{\sigma}} + \|\sqrt{\mathcal{T}_1}\| \right) \frac{R_0}{k+1},
\end{equation*}
and
\begin{equation*}\label{eq: complexity-bound-obj}
- \frac{1}{\sqrt{\sigma}} \|x^*\| \cdot \frac{R_0}{k+1}
\ \leq\ 
h(\by^{k+1}, \bw^{k+1}, \bz^{k+1}) 
\ \leq\ 
\left(3R_0 + \frac{1}{\sqrt{\sigma}} \|x^*\| \right) \frac{R_0}{k+1}.
\end{equation*}  
\end{proposition}

\subsection{A Dual HPR Method Incorporating the sGS Technique}
Note that the main computational bottleneck in Algorithm~\ref{alg:acc-pADMM-QP-0} lies in solving the subproblem with respect to the variables \( (y, w) \). To alleviate this difficulty, we apply the sGS technique~\citep{li2016schur,li2019block} to decouple the variables \( y \) and \( w \). Specifically, we define the self-adjoint positive semidefinite linear operator \( \mathcal{S}: \mathbb{R}^m \times \mathcal{W} \to \mathbb{R}^m \times \mathcal{W}\) as
\begin{equation}\label{def:S}
 \mathcal{S} = \sigma \begin{pmatrix}
\mathcal{S}_y & \mathbf{0} \\
\mathbf{0} & \mathcal{S}_w
\end{pmatrix},   
\end{equation}
where \( \mathcal{S}_y \) and \( \mathcal{S}_w \) are self-adjoint positive semidefinite linear operators such that \( \mathcal{S}_y + A A^{*} \) is positive definite on \( \mathbb{R}^m \). Then the sGS operator \( \widehat{\mathcal{S}}_{\rm sGS} : \mathbb{R}^m \times \mathcal{W} \to \mathbb{R}^m \times \mathcal{W} \) is defined as
\begin{equation}\label{def:sGS}
 \widehat{\mathcal{S}}_{\rm sGS} = 
\begin{pmatrix}
\mathbf{0} & -\sigma A Q \\
\mathbf{0} & \mathbf{0}
\end{pmatrix}
\begin{pmatrix}
\frac{1}{\sigma} (A A^{*} + \mathcal{S}_y)^{-1} & \mathbf{0} \\
\mathbf{0} & (\sigma Q^2 + Q + \sigma \mathcal{S}_w)^{-1}
\end{pmatrix}
\begin{pmatrix}
\mathbf{0} & \mathbf{0} \\
-\sigma Q A^{*} & \mathbf{0}
\end{pmatrix}
= 
\begin{pmatrix}
\widehat{\mathcal{S}}_{\rm sGS1} & \mathbf{0} \\
\mathbf{0} & \mathbf{0}
\end{pmatrix},   
\end{equation}
where \( \widehat{\mathcal{S}}_{\rm sGS1} := \sigma^2 A Q (\sigma Q^2 + Q + \sigma \mathcal{S}_w)^{-1} Q A^{*} \). According to \citep[Theorem 1]{li2019block}, the following proposition shows that incorporating the sGS operator enables an efficient block-wise update of \((\by^{k+1}, \bw^{k+1})\) for all $k \geq 0$.
\begin{proposition}\label{Prop:sGS}
Let \( \mathcal{S}_y \) and \( \mathcal{S}_w \) be two self-adjoint positive semidefinite linear operators on \( \mathbb{R}^m \) and \( \mathcal{W} \), respectively, such that \( \mathcal{S}_y + A A^{*} \) is positive definite. Let \( \mathcal{T}_1 = \mathcal{S} + \widehat{\mathcal{S}}_{\rm sGS} \), where \( \mathcal{S} \) and \( \widehat{\mathcal{S}}_{\rm sGS} \) are given in~\eqref{def:S} and~\eqref{def:sGS}. Then for any $k\geq 0$, Step 3 of Algorithm~\ref{alg:acc-pADMM-QP-0},
\[
(\by^{k+1}, \bw^{k+1}) = \underset{(y,w) \in \mathbb{R}^{m}\times \mathcal{W}}{\arg \min }\left\{L_\sigma(y,w,\bz^{k+1}; \bx^{k+1})+\frac{1}{2}\|(y,w)-(y^k,w^k)\|^2_{\cT_1}\right\}
\]
is equivalent to the following updates:
\[
\begin{cases}
\bw^{k+\frac{1}{2}} = \underset{w\in  \mathcal{W}}{\arg \min } \left\{ L_\sigma(y^k, w, \bz^{k+1}; \bx^{k+1}) + \frac{\sigma}{2} \|w - w^k\|^2_{\mathcal{S}_w} \right\}, \\
\by^{k+1} = \underset{y\in  \mathbb{R}^{m}}{\arg \min } \left\{ L_\sigma(y, \bw^{k+\frac{1}{2}}, \bz^{k+1}; \bx^{k+1}) + \frac{\sigma}{2} \|y - y^k\|^2_{\mathcal{S}_y} \right\}, \\
\bw^{k+1} = \underset{w\in  \mathcal{W}}{\arg \min } \left\{ L_\sigma(\by^{k+1}, w, \bz^{k+1}; \bx^{k+1}) + \frac{\sigma}{2} \|w - w^k\|^2_{\mathcal{S}_w} \right\}.
\end{cases}
\]
Moreover, \( \mathcal{T}_1 + \mathcal{H} \) is positive definite on \( \mathbb{R}^m \times \mathcal{W} \).
\end{proposition}
An HPR method incorporating the sGS technique is presented in Algorithm \ref{alg:acc-pADMM-QP}: 
\begin{algorithm}[H] 
	\caption{A dual HPR method for solving the restricted-Wolfe dual problem \eqref{QP-dual}}
	\label{alg:acc-pADMM-QP}
	\begin{algorithmic}[1] 
	\State {Input:  Let $\cS_y$ and  $\cS_w$ be two self-adjoint positive semidefinite linear operators on $\mathbb{R}^{m}$ and $\mathcal{W}$, respectively, such that $\cS_y+AA^{*}$ is positive definite. Denote $u=(y,w,z,x)$ and $\bu=(\by,\bw,\bz,\bx)$. Let \( u^0 = (y^0, w^0, z^0, x^0) \in \mathcal{U} \), and set  \( \sigma > 0 \).}
    		\For{\( k = 0, 1, 2, \ldots \)}
    \State {Step 1. $\displaystyle \bz^{k+1}=\underset{z \in \mathbb{R}^{n}}{\arg \min }\left\{L_\sigma(y^k, w^k,z; x^k)\right\}$;} 
    \State {Step 2. $\displaystyle  \bx^{k+1}={x}^k+\sigma (-Qw^{k}+A^{*}y^{k}+\bz^{k+1}-c)$;}
    \State {Step 3-1. $\displaystyle \bw^{k+\frac{1}{2}}=\underset{w \in \mathcal{W}}{\arg \min }\left\{L_\sigma(y^{k},w,\bz^{k+1}; \bx^{k+1})+\frac{\sigma}{2}\|w-w^k\|^2_{\cS_w}\right\}$;}
   \State {Step 3-2. $\displaystyle \by^{k+1}=\underset{y \in \mathbb{R}^{m}}{\arg \min }\left\{L_\sigma(y,\bw^{k+\frac{1}{2}},\bz^{k+1};\bx^{k+1})+\frac{\sigma}{2}\|y-y^k\|_{\cS_y}^2\right\}$;}	
   \State {Step 3-3. $\displaystyle \bw^{k+1}=\underset{w \in \mathcal{W}}{\arg \min }\left\{L_\sigma(\by^{k+1},w,  \bz^{k+1}; \bx^{k+1})+\frac{\sigma}{2}\|w-w^k\|^2_{\cS_w}\right\}$;}        
    \State {Step 4. $\displaystyle \hat{u}^{k+1}= 2 \bu^{k+1}-u^k$;}
	\State {Step 5. $\displaystyle u^{k+1}=\frac{1}{k+2}u^{0}+\frac{k+1}{k+2}\hat{u}^{k+1}$;}
    \EndFor
	\end{algorithmic}
\end{algorithm}

\begin{remark}
The proposed HPR method, as outlined in Algorithm \ref{alg:acc-pADMM-QP}, which incorporates the sGS technique, provides a unified and flexible framework for solving general CCQP problems of the form~\eqref{QP-primal}. In particular, when \(Q = \mathbf{0}\) and $\phi(\cdot)=\delta_{C}(\cdot)$, Algorithm~\ref{alg:acc-pADMM-QP} reduces exactly to the HPR method for LP introduced in~\citep{chen2024hpr}. Moreover, if the constraint \(Ax \in \mathcal{K}\) is absent—as in the case of the Lasso problem~\citep{li2018highly}—then Algorithm~\ref{alg:acc-pADMM-QP} simplifies to the HPR method with the update cycle \(z \to x \to w\). These special cases highlight the flexibility and broad applicability of the dual HPR method in Algorithm \ref{alg:acc-pADMM-QP}.
\end{remark}

\begin{remark}
In Algorithm~\ref{alg:acc-pADMM-QP}, the updates for \( \bw^{k+\frac{1}{2}} \) and \( \bw^{k+1} \) for $k\geq 0$ are restricted to the subspace \( \mathcal{W} = \operatorname{Range}(Q) \). Although it may seem more straightforward to update \( \bw^{k+\frac{1}{2}} \) and \( \bw^{k+1} \) in the full space \( \mathbb{R}^n \), doing so—particularly under a linearized ADMM framework—necessitates a proximal operator with a larger spectral norm, such as \( \mathcal{S}_w = \lambda_1(Q^2 + Q/\sigma) I_n - (Q^2 + Q/\sigma) \), to ensure convergence. A proximal operator with a large spectral norm typically results in slower convergence. By contrast, restricting the update to \( \mathcal{W} \) allows HPR-QP to employ a proximal operator with a smaller spectral norm, namely \( \mathcal{S}_w = Q(\lambda_1(Q) I_n - Q) \), which accelerates convergence while preserving theoretical guarantees.
\end{remark}

\subsection{An Easy-to-Implement Dual HPR Method}
While directly computing \( \bw^{k+\frac{1}{2}} \) and \( \bw^{k+1} \) for $k\geq0$ within the subspace \( \operatorname{Range}(Q) \) may appear computationally intensive, we show in this subsection that these updates can be performed efficiently without requiring explicit projection onto \( \operatorname{Range}(Q) \). For small‑scale problems or when \(Q\) has a favorable structure, one may set \(\mathcal{S}_w=0\) and solve the subproblems using direct solvers or preconditioned conjugate gradient methods \citep{liang2022qppal}. Here, for solving large‑scale general CCQP problems, we employ the proximal operator
\begin{equation}\label{def:S_w}
\mathcal{S}_w = Q(\lambda_Q I_n - Q),
\end{equation}
where \( \lambda_Q > 0 \) is a constant satisfying \( \lambda_Q \geq \lambda_1(Q) \).  Then, for any \( k \geq 0 \), the updates of \( \bw^{k+\frac{1}{2}} \) and \( \bw^{k+1} \) are given by
\begin{equation*}\label{update-w-true}
\left\{
\begin{aligned}
\displaystyle Q \bw^{k+\frac{1}{2}} &= \frac{1}{1 + \sigma \lambda_Q} Q \left( \sigma \lambda_Q w^{k} + \bx^{k+1} + \sigma(-Qw^{k} + A^* y^{k} + \bz^{k+1} - c) \right), \quad \bw^{k+\frac{1}{2}} \in \mathcal{W}, \\[1ex]
\displaystyle Q \bw^{k+1} &= \frac{1}{1 + \sigma \lambda_Q} Q \left( \sigma \lambda_Q w^{k} + \bx^{k+1} + \sigma(-Qw^{k} + A^* \by^{k+1} + \bz^{k+1} - c) \right), \quad \bw^{k+1} \in \mathcal{W}.
\end{aligned}
\right.
\end{equation*}
The following proposition—motivated by \citep[Proposition~4.1]{li2018qsdpnal}—demonstrates that this choice of \(\mathcal{S}_w\) enables efficient computation of the updates by constructing a shadow sequence, avoiding explicit projection onto \( \operatorname{Range}(Q) \).
\begin{proposition}\label{prop-Qw}
Let \( \mathcal{W} = \operatorname{Range}(Q)\) and $R\in \mathbb{R}^{n}$. To solve
\begin{equation}\label{eq:Qw-01}
Q \bw^{+} = \frac{1}{1 + \sigma \lambda_Q} \, Q R, \quad \bw^+ \in \mathcal{W},
\end{equation}
it suffices to set
\begin{equation}\label{eq:Qw-02}
\bw_{Q}^+ := \frac{1}{1 + \sigma \lambda_Q} \,  R.
\end{equation}
Then \( \bw^+ = \Pi_{\mathcal{W}}(\bw_{Q}^+) \) solves \eqref{eq:Qw-01}, and
\(Q \bw^{+}= Q \bw_{Q}^+, \, \langle \bw^+, Q \bw^+ \rangle = \langle \bw_{Q}^+, Q \bw_{Q}^+ \rangle.\)
\end{proposition}

Based on Proposition \ref{prop-Qw}, an easy-to-implement dual HPR method for solving the large-scale restricted-Wolfe dual problem \eqref{QP-dual} is presented in Algorithm \ref{alg:acc-pADMM-QP-im}.
\begin{algorithm}[H] 
	\caption{An easy-to-implement  dual HPR method for the restricted-Wolfe dual problem \eqref{QP-dual}}
	\label{alg:acc-pADMM-QP-im}
	\begin{algorithmic}[1] 
	\State {Input: Let \(\mathcal{S}_w\) be defined as in \eqref{def:S_w}, and let \(\mathcal{S}_y\) be a self-adjoint positive semidefinite linear operator on \(\mathbb{R}^m\) such that \(\mathcal{S}_y + A A^*\) is positive definite. Denote $u_Q=(y,w_Q,z,x)$ and $\bu_{Q}=(\by,\bw_{Q},\bz,\bx)$. Let \( u_{Q}^0 = (y^0, w_{Q}^0, z^0, x^0) \in \mathbb{R}^{m}\times \mathbb{R}^{n}\times\mathbb{R}^{n} \times\mathbb{R}^{n} \), and set \( \sigma > 0 \).}
    		\For{\( k = 0, 1, 2, \ldots \)}
    \State {Step 1. $\displaystyle \bz^{k+1}=\underset{z \in \mathbb{R}^{n}}{\arg \min }\left\{L_\sigma(y^k, w^k_{Q},z; x^k)\right\} $;} 
    \State {Step 2. $\displaystyle  \bx^{k+1}={x}^k+\sigma (-Qw_{Q}^{k}+A^{*}y^{k}+\bz^{k+1}-c) $;}\vspace{3pt}
    \State {Step 3-1. $\displaystyle \bw_{Q}^{k+\frac{1}{2}}=\frac{1}{1 + \sigma \lambda_Q}  \left( \sigma\lambda_Q w_{Q}^{k} + \bx^{k+1} + \sigma(-Qw^{k}_{Q} + A^* y^{k} + \bz^{k+1} - c) \right)$;}\vspace{3pt}
   \State {Step 3-2. $\displaystyle \by^{k+1}=\underset{y \in \mathbb{R}^{m}}{\arg \min }\left\{L_\sigma(y,\bw_{Q}^{k+\frac{1}{2}},\bz^{k+1};\bx^{k+1})+\frac{\sigma}{2}\|y-y^k\|_{\cS_y}^2\right\}$;}	\vspace{3pt}
   \State {Step 3-3. $\displaystyle \bw_{Q}^{k+1}=\frac{1}{1 + \sigma \lambda_Q}  \left( \sigma \lambda_Q w_{Q}^{k} + \bx^{k+1} + \sigma(-Qw^{k}_{Q} + A^* \by^{k+1} + \bz^{k+1} - c) \right)$;} \vspace{3pt}       
    \State {Step 4. $\displaystyle \hat{u}_{Q}^{k+1}= 2 \bu_{Q}^{k+1}-u_{Q}^k$;}\vspace{3pt}
    \State {Step 5. $\displaystyle u_{Q}^{k+1}=\frac{1}{k+2}u_{Q}^{0}+\frac{k+1}{k+2}\hat{u}_{Q}^{k+1}$;}
    \EndFor
	\end{algorithmic}
\end{algorithm}

\begin{theorem}\label{th:convergence-HPR-imple}
Suppose Assumption~\ref{ass: CQ-QP} holds. Then the sequence \(\{(\by^k,\,Q\bw_Q^k,\,\bz^k,\,\bx^k)\}\) generated by Algorithm~\ref{alg:acc-pADMM-QP-im} is equivalent to the sequence \(\{(\by^k,\,Q\bw^k,\,\bz^k,\,\bx^k)\}\) produced by Algorithm~\ref{alg:acc-pADMM-QP} starting from the same initial point. Moreover, both sequences converge to the point \((y^*,\,Qw^*,\,z^*,\,x^*)\), where \((y^*,\,w^*,\,z^*)\) solves problem~\eqref{QP-dual} and \(x^*\) solves problem~\eqref{QP-primal}.
\end{theorem}

\begin{proof}
According to the optimality conditions of the subproblems in Algorithm~\ref{alg:acc-pADMM-QP}, for any \(k \geq 0\), we have:
\begin{equation}\label{update}
\begin{cases}
\mathbf{0} \in -\partial \phi^{*}(-\bz^{k+1}) + x^k + \sigma(-Qw^k + A^{*}y^k + \bz^{k+1} - c), \\
\bx^{k+1} = x^k + \sigma(-Qw^k + A^* y^k + \bz^{k+1} - c), \\
Q \bw^{k+\frac{1}{2}} = \frac{1}{1 + \sigma \lambda_Q} Q \left( \sigma \lambda_Q w^k + \bx^{k+1} + \sigma(-Qw^k + A^* y^k + \bz^{k+1} - c) \right),\quad \bw^{k+\frac{1}{2}} \in \mathcal{W}, \\
\mathbf{0} \in -\partial \delta^{*}_{\mathcal{K}}(-\by^{k+1}) + A\bx^{k+1} + \sigma A(-Q\bw^{k+\frac{1}{2}} + A^*\by^{k+1} + \bz^{k+1} - c) + \sigma \cS_y(\by^{k+1} - y^k), \\
Q \bw^{k+1} = \frac{1}{1 + \sigma \lambda_Q} Q \left( \sigma \lambda_Q w^k + \bx^{k+1} + \sigma(-Qw^k + A^* \by^{k+1} + \bz^{k+1} - c) \right),\quad \bw^{k+1} \in \mathcal{W}.
\end{cases}
\end{equation}
Based on \eqref{update} and Proposition~\ref{prop-Qw}, we can derive that the sequence \(\{(\by^k,\,Q\bw_Q^k,\,\bz^k,\,\bx^k)\}\) generated by Algorithm~\ref{alg:acc-pADMM-QP-im} is equivalent to \(\{(\by^k,\,Q\bw^k,\,\bz^k,\,\bx^k)\}\) from Algorithm~\ref{alg:acc-pADMM-QP}, under the same initialization.

Furthermore, Proposition~\ref{Prop:sGS} establishes that Algorithm~\ref{alg:acc-pADMM-QP} is a special case of Algorithm~\ref{alg:acc-pADMM-QP-0} with \(\cT_1 = \cS +  \widehat{\mathcal{S}}_{\rm sGS}\). Thus, the convergence result in Proposition~\ref{Prop:convergence} applies, and the sequence \(\{(\by^k,\,\bw^k,\,\bz^k,\,\bx^k)\}\) from Algorithm~\ref{alg:acc-pADMM-QP} converges to the optimal solution of the primal-dual pair \eqref{QP-primal} and \eqref{QP-dual}. This completes the proof.

\end{proof}

By substituting \(\mathcal{T}_1 = \mathcal{S} + \widehat{\mathcal{S}}_{\rm sGS}\), where \(\mathcal{S}\) and \(\widehat{\mathcal{S}}_{\rm sGS}\) are defined in \eqref{def:S} and \eqref{def:sGS}, respectively, into the definition of \(\mathcal{M}\) in \eqref{def:M-QP}, we obtain the explicit form of \(\mathcal{M}\) used in the dual HPR method (Algorithm~\ref{alg:acc-pADMM-QP-im}):
\begin{equation}
\label{eq:metric-M}
\mathcal{M} = 
\displaystyle\begin{bmatrix}
\sigma A_{Q}^{*}A_{Q}+\cS+ \widehat{\mathcal{S}}_{\rm sGS}  & \quad  \mathbf{0} & \quad A^{*}_{Q} \\
\mathbf{0} & \quad \mathbf{0} &  \quad   \mathbf{0} \\
A_{Q} & \quad \mathbf{0} &    \quad \frac{1}{\sigma} I_n\\
\end{bmatrix}.
\end{equation}
Combining Propositions~\ref{prop:complexity-acc-pADMM}, \ref{Prop:sGS}, and Theorem~\ref{th:convergence-HPR-imple}, we derive the following complexity result for Algorithm~\ref{alg:acc-pADMM-QP-im}:

\begin{theorem}\label{th:complexity-HPR-QP}
Suppose Assumption~\ref{ass: CQ-QP} holds. Let \(\{\bu_{Q}^k\} = \{(y^k,\,w_{Q}^k,\,z^k,\,x^k)\}\) and \(\{u_{Q}^k\} = \{(y^k,\,w_{Q}^k,\,z^k,\,x^k)\}\) be the sequences generated by Algorithm~\ref{alg:acc-pADMM-QP-im}, and let \( u^* = (y^*, w^*, z^*, x^*) \) be a solution to the KKT system~\eqref{eq:KKT-QP}. Then, for all \(k \geq 0\), the following complexity bounds hold with \(R_0 = \|u_{Q}^0 - u^*\|_{\mathcal{M}}\):
\begin{equation*}
    \begin{aligned}
\|\bu_{Q}^{k+1} - u_{Q}^k\|_{\mathcal{M}} &\le \frac{R_0}{k+1}, \label{eq:complexity-M}\\
\|\mathcal{R}(\bar u_{Q}^{k+1})\| &\le \left(\frac{\sigma\|A^{*}_{Q}\| + 1}{\sqrt{\sigma}} + \|\sqrt{\cS+ \widehat{\mathcal{S}}_{\rm sGS}}\|\right)\frac{R_0}{k+1}, 
\\
- \frac{\|x^*\|}{\sqrt{\sigma}} \cdot \frac{R_0}{k+1}
&\le h(\by^{k+1}, \bw_{Q}^{k+1}, \bz^{k+1})
\le \left(3R_0 + \frac{\|x^*\|}{\sqrt{\sigma}}\right)\frac{R_0}{k+1}.
\end{aligned}
\end{equation*}
\end{theorem}

\section{HPR-QP: A Dual HPR Method for CCQP}\label{sec:3}
In this section, we present HPR-QP, as outlined in Algorithm~\ref{alg:HPR_QP}, for solving problem~\eqref{QP-dual}, incorporating an adaptive restart strategy and a dynamic update rule for \( \sigma \).

\begin{algorithm}[ht!]
\caption{HPR-QP: A dual HPR method for the  CCQP \eqref{QP-dual}}
\label{alg:HPR_QP}
\begin{algorithmic}[1]
\State \textbf{Input:} Let \(\mathcal{S}_w\) be defined as in \eqref{def:S_w}, and let \(\mathcal{S}_y\) be a self-adjoint, positive semidefinite operator on \(\mathbb{R}^m\) such that \(\mathcal{S}_y + A A^*\) is positive definite. Let $u_{Q} = (y, w_{Q}, z, x)$, $\bu_{Q} = (\by, \bw_{Q}, \bz, \bx)$, and initial point $u_{Q}^{0,0} = (y^{0,0}, w_{Q}^{0,0}, z^{0,0}, x^{0,0}) \in \mathbb{R}^{m}\times \mathbb{R}^{n}\times\mathbb{R}^{n} \times\mathbb{R}^{n}$.
\State \textbf{Initialization:} Set outer loop counter $r = 0$, total iteration counter $k = 0$, and initial penalty parameter $\sigma_0 > 0$.
\vspace{3pt}

\Repeat
    \State Initialize inner loop: set inner counter $t = 0$;
    \Repeat
        \State $\displaystyle \bz^{r,t+1} = \arg\min_{z \in \mathbb{R}^n} L_{\sigma_r}(y^{r,t}, w^{r,t}_{Q}, z; x^{r,t})$;
        \State $\displaystyle \bx^{r,t+1} = x^{r,t} + \sigma_r (-Q w_{Q}^{r,t} + A^* y^{r,t} + \bz^{r,t+1} - c)$;
        \State $\displaystyle \bw_{Q}^{r,t+\frac{1}{2}} = \frac{1}{1 + \sigma_r \lambda_Q} \left( \sigma_r \lambda_Q w_{Q}^{r,t} + \bx^{r,t+1} + \sigma_r(-Qw_{Q}^{r,t} + A^* y^{r,t} + \bz^{r,t+1} - c) \right)$;
        \State $\displaystyle \by^{r,t+1} = \arg\min_{y \in \mathbb{R}^m} \left\{ L_{\sigma_r}(y, \bw_{Q}^{r,t+\frac{1}{2}}, \bz^{r,t+1}; \bx^{r,t+1}) + \frac{\sigma_r}{2} \|y - y^{r,t}\|_{\cS_y}^2 \right\}$;
        \State $\displaystyle \bw_{Q}^{r,t+1} = \frac{1}{1 + \sigma_r \lambda_Q} \left( \sigma_r \lambda_Q w_{Q}^{r,t} + \bx^{r,t+1} + \sigma_r(-Qw_{Q}^{r,t} + A^* \by^{r,t+1} + \bz^{r,t+1} - c) \right)$;
        \State $\displaystyle \hat{u}_{Q}^{r,t+1} = 2\bu_{Q}^{r,t+1} - u_{Q}^{r,t}$;
        \State $\displaystyle u_{Q}^{r,t+1} = \frac{1}{k+2} u_{Q}^{r,0} + \frac{k+1}{k+2} \hat{u}_{Q}^{r,t+1}$;
        \State $t \gets t + 1$, $k \gets k + 1$;
    \Until restart or termination criteria are met;
    \State \textbf{Restart:} Set $\tau_r = t$, $u_{Q}^{r+1,0} = \bu_{Q}^{r,\tau_r}$;
    \State $\displaystyle \sigma_{r+1} = \textbf{SigmaUpdate}(\bu_{Q}^{r,\tau_r}, u_{Q}^{r,0}, \cS_y, \cS_w, A, Q)$;
    \State $r \gets r + 1$;
\Until termination criteria are met;
\State \textbf{Output:} $\{\bu^{r,t}\}$.
\end{algorithmic}
\end{algorithm}

\subsection{Efficient Solution of Subproblems}
We first detail the update formulas for each subproblem in HPR-QP (Steps 6–10). Specifically, for any \( r \geq 0 \) and \( t \geq 0 \), the update of \( \bz^{r,t+1} \) is given by
\begin{equation}\label{update-z}
    \bz^{r,t+1} = \frac{1}{\sigma_r} \left( \operatorname{Prox}_{\sigma_r\phi}(r_{z}^{r,t}) - r_{z}^{r,t} \right), 
\end{equation}
where \( r_{z}^{r,t} = x^{r,t} + \sigma_r (-Qw_{Q}^{r,t} + A^* y^{r,t} - c) \). The corresponding update for \( \bx^{r,t+1} \) becomes
\begin{equation}\label{update-x}
\bx^{r,t+1} = x^{r,t} + \sigma_r(-Qw_{Q}^{r,t} + A^* y^{r,t} + \bz^{r,t+1} - c) = \operatorname{Prox}_{\sigma_{r}\phi}(r_{z}^{r,t}).
\end{equation}
Furthermore, the updates for \( \bw_Q^{r,t+\frac{1}{2}} \) and  \( \bw_Q^{r,t+1} \)  can be simplified as follows
\begin{equation}\label{update-w}
\left\{
\begin{aligned}
\bw_{Q}^{r,t+\frac{1}{2}} &= \frac{1}{1 + \sigma_r \lambda_Q}  \left( \sigma_r \lambda_Q w_{Q}^{r,t} + \bx^{r,t+1} + \sigma_r(-Qw_{Q}^{r,t} + A^* y^{r,t} + \bz^{r,t+1} - c) \right) \\[1ex]
&= \frac{1}{1 + \sigma_r \lambda_Q} \left( \sigma_r \lambda_Q w_{Q}^{r,t} + 2\bx^{r,t+1} - x^{r,t} \right) = \frac{1}{1 + \sigma_r \lambda_Q} \left( \sigma_r \lambda_Q w_{Q}^{r,t} + \hx^{r,t+1} \right), \\[1.5ex]
 \bw_{Q}^{r,t+1} &=\frac{1}{1 + \sigma_r \lambda_Q} \left( \sigma_r \lambda_Q w_{Q}^{r,t} + \bx^{r,t+1} + \sigma_r(-Qw_{Q}^{r,t} + A^* \by^{r,t+1} + \bz^{r,t+1} - c) \right)\\
&=\left( \bw_{Q}^{r,t+\frac{1}{2}} + \frac{\sigma_r}{1 + \sigma_r \lambda_Q} A^*(\by^{r,t+1} - y^{r,t}) \right).
\end{aligned}
\right.
\end{equation}
To simplify the solution of the subproblem with respect to \( y \), we choose the proximal operator
\begin{equation}\label{def:S_y}
\mathcal{S}_y = \lambda_A I_m - A A^*,
\end{equation}
where  \( \lambda_A > 0 \) is a constant such that  \( \lambda_A \geq \lambda_{1}(A A^*) \). 
 The updates for  \( \by^{r,t+1} \) are then given by
\begin{equation}\label{update-y}
\begin{aligned}
\by^{r,t+1} &= \frac{1}{\sigma_r \lambda_A } \left( \Pi_{\cK}(r_y^{r,t}) - r_y^{r,t} \right), \\[1.5ex]
\end{aligned}
\end{equation}
where 
\[
\begin{aligned}
r_{y}^{r,t} &= A\left( \bx^{r,t+1} + \sigma_r (-Q \bw_{Q}^{r,t+\frac{1}{2}} + A^* y^{r,t} + \bz^{r,t+1} - c) \right) - \sigma_r \lambda_A y^{r,t} \\
&= A(\hx^{r,t+1} + \sigma_r(Qw_{Q}^{r,t} - Q\bw_{Q}^{r,t+\frac{1}{2}})) - \sigma_r \lambda_A y^{r,t}.
\end{aligned}
\]
Moreover, by combining the update formulas in \eqref{update-x}, \eqref{update-w},  and \eqref{update-y}, we find that it is not necessary to compute \( \bz^{r,t+1} \) at every iteration. Instead, \( \bz^{r,t+1} \) needs to be evaluated via \eqref{update-z} only when checking the stopping criteria. This further improves computational efficiency without affecting the correctness of the method.

Finally, to fully leverage the parallel computing capabilities of GPUs, we implement custom CUDA kernels for the updates in \eqref{update-z}, \eqref{update-x}, \eqref{update-w}, and \eqref{update-y}. For matrix-vector multiplications, we selectively use either custom kernel implementations or the \texttt{cusparseSpMV()} routine from the cuSPARSE library, which applies the \texttt{CUSPARSE\_SPMV\_CSR\_ALG2} algorithm. Our custom implementations reduce the number of kernel launches, which is beneficial for performance, especially when minimizing launch overhead. However, they may be less efficient than cuSPARSE when handling dense or moderately sparse matrices. Therefore, when A and Q are not highly sparse, we prefer \texttt{cusparseSpMV()} to ensure optimal performance.

\subsection{An Adaptive Restart Strategy}
Inspired by the restart strategy proposed in PDLP~\citep{applegate2021practical,lu2023cupdlp,lu2023cupdlp_c}, HPR-LP~\citep{chen2024hpr} introduced an adaptive restart mechanism based on the \( O(1/k) \) iteration complexity of the HPR method, which has demonstrated practical success on large-scale LP problems. Motivated by this, we extend this adaptive restart strategy to the CCQP~\eqref{QP-primal} by defining a suitable merit function grounded in the complexity result established in Theorem~\ref{th:complexity-HPR-QP}.
Specifically, we define the following idealized merit function:
\[
R_{r,t} := \|u_{Q}^{r,t} - u^*\|_{\mathcal{M}} \quad \forall r \geq 0, \ t \geq 0,
\]
where \( u^* \) denotes any solution to the KKT system~\eqref{eq:KKT-QP}. Note that \( R_{r,0} \) corresponds to the upper bound given by the complexity result in Theorem~\ref{th:complexity-HPR-QP} at the start of the \( r \)-th outer iteration. Since \( u^* \) is unknown in practice, we approximate \( R_{r,t} \) using the computable surrogate:
\[
\widetilde{R}_{r,t} := \|u_{Q}^{r,t} - \bar{u}_{Q}^{r,t+1}\|_{\mathcal{M}}.
\]
Based on this surrogate merit function, we introduce the following adaptive restart criteria for the HPR-QP method:
\begin{enumerate}
    \item {Sufficient decay:}
    \begin{equation}\label{eq:restart_1}
        \widetilde{R}_{r,t+1} \leq \alpha_1 \widetilde{R}_{r,0};
    \end{equation}
    
    \item {Insufficient local progress despite overall decay:}
    \begin{equation}\label{eq:restart_2}
        \widetilde{R}_{r,t+1} \leq \alpha_2 \widetilde{R}_{r,0} \quad \text{and} \quad \widetilde{R}_{r,t+1} > \widetilde{R}_{r,t};
    \end{equation}

    \item {Excessively long inner loop:}
    \begin{equation}\label{eq:restart_3}
        t \geq \alpha_3 k;
    \end{equation}
\end{enumerate}
where \( \alpha_1 \in (0, \alpha_2) \), \( \alpha_2 \in (0,1) \), and \( \alpha_3 \in (0,1) \) are user-defined parameters. Whenever any of the above conditions is satisfied, we terminate the current inner loop and begin a new outer iteration by setting \( u^{r+1,0} = \bar{u}^{r,\tau_r} \) and updating the penalty parameter \( \sigma_{r+1} \) accordingly.

\subsection{An Update Strategy for \texorpdfstring{$\sigma$}{Sigma}}
Now, we describe the update rule for the penalty parameter \( \sigma \) in HPR-QP. This design is also motivated by the complexity results of the HPR method with semi-proximal terms, as established in Theorem~\ref{th:complexity-HPR-QP}. Specifically, at the beginning of the \( (r+1) \)-th outer iteration, we determine \( \sigma_{r+1} \) by solving the following minimization problem:
\begin{equation}\label{eq:sigma-0}
\sigma_{r+1} = \arg \min_{\sigma > 0} \| u_{Q}^{r+1,0} - u^* \|_{\mathcal{M}}^2,
\end{equation}
where \( u^* \) denotes a solution to the KKT system~\eqref{eq:KKT-QP}, and the metric \( \| u_{Q}^{r+1,0} - u^* \|_{\mathcal{M}} \) corresponds to the upper bound derived in Theorem~\ref{th:complexity-HPR-QP}. The rationale behind this formulation is that minimizing this bound is expected to yield a smaller quantity \( \|u_{Q}^{r+1,t} - \bar{u}_{Q}^{r+1,t+1}\|_{\mathcal{M}} \) for all \( t \geq 0 \), which, in turn, reduces the KKT residual \( \|\mathcal{R}(\bar{u}_{Q}^{r+1,t+1})\| \). This insight guides the dynamic adjustment of \( \sigma \) to enhance the practical performance of the method.

Substituting \eqref{eq:metric-M} and \eqref{def:S_y} into \eqref{eq:sigma-0}, the update rule for \( \sigma_{r+1} \) reduces to solving the following one-dimensional minimization problem:
\begin{equation}\label{eq:update-sigma}
\sigma_{r+1} = \arg \min_{\sigma > 0} \left\{ f (\sigma):=\theta_1 \sigma + \frac{\theta_2}{\sigma} + \frac{\sigma^2 \theta_3}{1 + \lambda_Q \sigma} \right\},
\end{equation}
where the coefficients \( \theta_1, \theta_2, \theta_3 \) are given by 
\[
\begin{aligned}
\theta_1 &= \lambda_{A} \|y^{r+1,0} - y^*\|^2 + \lambda_Q \|w_{Q}^{r+1,0} - w^*\|_Q^2 - 2 \langle Q(w_{Q}^{r+1,0} - w^*), A^*(y^{r+1,0} - y^*) \rangle, \\
\theta_2 &= \|x^{r+1,0} - x^*\|^2, \quad  
\theta_3 = \|A^* y^{r+1,0} - A^* y^*\|_Q^2. 
\end{aligned}
\]
When coefficients \( \theta_i > 0, \,i=2,3 \), the function \( f(\sigma) \) is strictly convex on \( \sigma > 0 \), since its second derivative satisfies
\[
f^{\prime\prime}(\sigma) = \frac{2 \theta_2}{\sigma^3} + \frac{2 \theta_3}{(1 + \lambda_Q \sigma)^3} > 0 \quad \text{for all } \sigma > 0.
\]
This guarantees both the existence and uniqueness of the optimal solution to problem \eqref{eq:update-sigma}. This scalar optimization problem can be efficiently solved using the golden section search method \citep{gerald2004applied}. Moreover, since the exact values of \( \theta_i \) involve the unknown solution \( u^* \), we approximate them in practice as follows:
\begin{equation}\label{eq:theta}
\begin{aligned}
\widetilde{\theta}_1 &= \lambda_{A} \|\by^{r,\tau_r} - y^{r,0}\|^2 + \lambda_Q \|\bw_{Q}^{r,\tau_r} - w_{Q}^{r,0}\|_Q^2 - 2 \langle Q(\bw_{Q}^{r,\tau_r} - w_{Q}^{r,0}), A^*(\by^{r,\tau_r} - y^{r,0}) \rangle, \\
\widetilde{\theta}_2 &= \|\bx^{r,\tau_r} - x^{r,0}\|^2, \quad 
\widetilde{\theta}_3 = \|A^* \by^{r,\tau_r} - A^* y^{r,0}\|_Q^2.
\end{aligned}
\end{equation}
Using these approximations, the updated penalty parameter \( \sigma \) is obtained by solving the following scalar optimization problem:
\begin{equation}\label{eq:update-sigma-approx}
\sigma_{\rm new} = \arg \min_{\sigma > 0} \left\{ f (\sigma)=\widetilde{\theta}_1 \sigma + \frac{\widetilde{\theta}_2}{\sigma} + \frac{\sigma^2 \widetilde{\theta}_3}{1 + \lambda_Q \sigma} \right\}.
\end{equation}
To further stabilize the update and prevent the approximations from deviating significantly from the true values, we apply an exponential smoothing scheme. The complete update procedure for \( \sigma \) is summarized in Algorithm~\ref{alg:sigma}:
\begin{algorithm}[H]
\caption{\textbf{SigmaUpdate}}
\label{alg:sigma}
\begin{algorithmic}[1]
    \State \textbf{Input:} \( \bu_{Q}^{r,\tau_r}, u_{Q}^{r,0}, \mathcal{S}_y, \mathcal{S}_w, A, Q \).
    \State Compute \( \widetilde{\theta}_1, \widetilde{\theta}_2, \widetilde{\theta}_3 \) as defined in~\eqref{eq:theta};
        \State Ensure numerical stability:
    \(
        \widetilde{\theta}_i \leftarrow \max\bigl(\widetilde{\theta}_i, 10^{-12}\bigr),\; i=1,2
    \);
    \State Solve problem \eqref{eq:update-sigma-approx} using golden section search method to obtain \( \sigma_{\mathrm{new}} \);
      \State Compute smoothing factor: \( \beta = \exp\left(- {\widetilde{R}_{r,\tau_r-1}}/{\widetilde{R}_{0,\tau_0-1}} \right) \);
    \State \textbf{Output:} \( \sigma_{r+1} = \exp\big( \beta \log(\sigma_{\mathrm{new}}) + (1 - \beta) \log(\sigma_r) \big) \).
\end{algorithmic}
\end{algorithm}

\section{Numerical Experiments}\label{sec:4}
In this section, we evaluate the performance of the Julia implementation of HPR-QP on a GPU and compare it against several state-of-the-art CQP solvers, including PDQP \citep{lu2023practical}, SCS~\citep{o2016conic,o2021operator}, CuClarabel~\citep{chen2024cuclarabel}, and Gurobi~\citep{gurobi}. The details of the experimental setup are provided in Section~\ref{Exper-setup}. We report numerical results on a diverse set of benchmark problems: Section~\ref{numerical-MM} presents results on the classical Maros-M\'esz\'aros data set~\citep{maros1999repository}; Section~\ref{numerical-Sy-QP} covers large-scale synthetic CQP instances from six problem classes introduced in~\citep{stellato2020osqp}; Section~\ref{numerical-Lasso} focuses on Lasso regression problems; and Section~\ref{numerical-QAP} reports results on large-scale CQP problems arising from relaxations of QAPs.

\subsection{Experimental Setup}\label{Exper-setup}

\paragraph{Solvers and Computing Environment.}
The {HPR-QP} solver is implemented in Julia~\citep{bezanson2017julia} with GPU acceleration enabled via CUDA. For comparison, both PDQP\footnote{https://github.com/jinwen-yang/PDQP.jl, downloaded in April 2025} and CuClarabel\footnote{https://github.com/oxfordcontrol/Clarabel.jl/tree/CuClarabel, version 0.10.0} are also implemented in Julia with support for CUDA-based GPU execution. In contrast, SCS\footnote{https://github.com/jump-dev/SCS.jl, version 2.1.0} is developed in C/C++ with a Julia interface and uses GPU acceleration through its indirect solver, where all matrix operations are performed on the GPU. Gurobi (version 12.0.2, academic license) is executed on the CPU. All solvers are benchmarked on a SuperServer SYS-420GP-TNR equipped with an NVIDIA A100-SXM4-80GB GPU, an Intel Xeon Platinum 8338C CPU @ 2.60 GHz, and 256 GB of RAM. The experiments are conducted on Ubuntu 24.04.2 LTS.

\paragraph{Preconditioning.}
Similar to PDQP~\citep{lu2023practical} and HPR‑LP \citep{chen2024hpr}, {HPR‑QP} also employs a diagonal preconditioning strategy when \(Q\) is available in explicit matrix form. In this setting, we perform 10 iterations of Ruiz equilibration~\citep{ruiz2001scaling}, followed by the diagonal preconditioning method of Pock and Chambolle~\citep{pock2011diagonal} with parameter \(\alpha = 1\). If \(Q\) is provided implicitly (i.e., as a matrix-free operator), preconditioning is not applied in the current implementation.

\paragraph{Parameter Setting.}
After preconditioning, we estimate \(\lambda_A\) and \(\lambda_Q\) using the power method~\citep{golub2013matrix}. {HPR-QP} is initialized at the origin, and the initial penalty parameter is set to \(\sigma_0 = {\|b\|}/{\|c\|},\) when both \(\|b\|\) and \(\|c\|\) lie within the range \([10^{-16}, 10^{16}]\); here, \(b := \max(|l|, |u|)\) is taken componentwise, treating any infinite entries in \(l\) or \(u\) as zero when computing the norm. Otherwise, we set \(\sigma_0 = 1\) to ensure numerical stability. The adaptive restart mechanism follows the criteria in \eqref{eq:restart_1}--\eqref{eq:restart_3}, with parameters \(\alpha_1 = 0.2\), \(\alpha_2 = 0.8\), and \(\alpha_3 = 0.5\). If the residual ratio satisfies \(\widetilde{R}_{r,\tau_r-1} / \widetilde{R}_{0,\tau_0-1} \leq 0.1\), then \(\alpha_3\) is tightened to 0.2. The penalty parameter \(\sigma\) is dynamically updated according to Algorithm~\ref{alg:sigma}. All other solvers are executed using their default parameter settings.

\paragraph{Termination and Time Limit.}
We terminate HPR-QP when the stopping criteria, similar to those in PDQP~\citep{lu2023practical}, are satisfied for a given tolerance \(\varepsilon \in (0, \infty)\). Specifically, the solver stops when the following conditions hold:
\[
\eta_{\mathrm{gap}} = 
\frac{\Big| \frac{1}{2}\langle x, Q x \rangle + \langle c, x \rangle + \phi(x) - \big( -\frac{1}{2}\langle x, Q x \rangle - \delta_{\mathcal{K}}^*(-y) - \phi^*(-z) \big) \Big|}
{1 + \max\Big( \big| \frac{1}{2}\langle x, Q x \rangle + \langle c, x \rangle + \phi(x) \big|, \big| \frac{1}{2}\langle x, Q x \rangle + \delta_{\mathcal{K}}^*(-y) + \phi^*(-z) \big| \Big)} \leq \varepsilon,
\]
\[
\eta_{\mathrm{p}} = \frac{\|Ax - \Pi_{\mathcal{K}}(Ax)\|_{\infty}}{1 + \max\left( \|b\|_{\infty}, \|Ax\|_{\infty} \right)} \leq \varepsilon,
\quad
\eta_{\mathrm{d}} = \frac{\| -Qx + A^* y + z - c \|_{\infty}}{1 + \max\left( \|c\|_{\infty}, \|A^* y\|_{\infty}, \|Qx\|_{\infty} \right)} \leq \varepsilon.
\]
All other solvers are run with their respective default stopping conditions. We evaluate the performance of each solver using three accuracy levels: \(\varepsilon = 10^{-4}\), \(\varepsilon = 10^{-6}\), and \(\varepsilon = 10^{-8}\). The numerical results for $\varepsilon=10^{-4}$ are reported in Appendix \ref{Append-numerical-result}. Additionally, a time limit of 3600 seconds is imposed on all algorithms for each problem instance across all data sets, except for the extremely large-scale QAP relaxations.

\paragraph{Shifted Geometric Mean.}
To evaluate solver performance over multiple problems, we use the shifted geometric mean (SGM), as in Mittelmann's benchmarks. For a shift \(\Delta = 10\), the SGM10 is defined as
\(
\left( \prod_{i=1}^n (t_i + \Delta) \right)^{1/n} - \Delta,
\)
where \(t_i\) denotes the solve time (or iteration count) in seconds for the \(i\)-th instance. Unsolved instances are assigned a time limit. Following PDQP~\citep{lu2023cupdlp}, we exclude the time spent on data loading and preconditioning from the measured solve time. In HPR‑QP, however, the time consumed by the power method is counted as part of the solve time, while PDQP does not include it. For SCS~\citep{o2016conic,o2021operator} and CuClarabel~\citep{chen2024cuclarabel}, we record only the solve time and exclude any setup overhead.  For Gurobi~\citep{gurobi}, we report the solve time as provided by the solver itself. 

\paragraph{Absolute Performance Profile.}
We also evaluate solver performance using the absolute performance profile \( f_s^a: \mathbb{R}_{+} \rightarrow [0,1] \), which represents the fraction of problems solved by solver \( s \) within time \( \tau \). It is defined as
\[
f_s^a(\tau) = \frac{1}{N} \sum_{p} \mathcal{I}_{\leq \tau}(t_{p, s}),
\]
where \( \mathcal{I}_{\leq \tau}(u) = 1 \) if \( u \leq \tau \), and \( 0 \) otherwise.

\subsection{Numerical Results on the Maros-M\'esz\'aros Data Set}\label{numerical-MM}
We begin by evaluating all tested algorithms on the Maros-M\'esz\'aros data set~\citep{maros1999repository}, a standard benchmark for CQP problems comprising 137 instances.\footnote{The instance \texttt{VALUES} is excluded from our evaluation, as Gurobi reports it to be non-convex.} Table~\ref{tab:d-p1-p2} demonstrates the superior performance of {HPR-QP} over the HPR method applied to the primal reformulations~\eqref{QP-primal-s} and~\eqref{QP-primal-sv}, even when both use similar adaptive restart and penalty update strategies. Notably, our primal HPR variants also outperform PDQP in terms of SGM10 for both runtime and iteration count, further underscoring the efficiency and robustness of our dual HPR framework.

\begin{table}[H]
  \centering
  \caption{Performance of HPR methods and PDQP on 15 Maros-M\'esz\'aros instances with tolerance \(10^{-8}\). {HPR} (p1) and {HPR} (p2) correspond to CCQP formulations~\eqref{QP-primal-s} and~\eqref{QP-primal-sv}, respectively.}
  \label{tab:d-p1-p2}
  \renewcommand{\arraystretch}{1.2}
  \begin{tabular}{l *{4}{rr}}
    \toprule
    \multirow{2}{*}{Instance}
      & \multicolumn{2}{c}{HPR-QP}
      & \multicolumn{2}{c}{HPR (p1)}
      & \multicolumn{2}{c}{HPR (p2)}
      & \multicolumn{2}{c}{PDQP} \\
    \cmidrule(lr){2-3}
    \cmidrule(lr){4-5}
    \cmidrule(lr){6-7}
    \cmidrule(lr){8-9}
      & Iter  & Time    & Iter     & Time    & Iter     & Time    & Iter      & Time    \\
    \midrule
    CVXQP3\_L   &  114\,400 &   5.4   &  110\,200  &   3.4   &   244\,200 &   9.2   &   136\,512 &   37.4  \\
    DUALC1      &  1\,900   &   0.1   &   6\,800   &   0.7   &   7\,200   &   0.6   &  10\,560   &   3.5   \\
    Q25FV47     &215\,100   &  32.4   & 379\,400   &  37.8   & 531\,600   &  49.9   & 891\,552   & 245.6   \\
    QBANDM      & 23\,500   &   1.4   & 136\,600   &   5.9   & 213\,900   &   9.9   & 139\,968   &  37.6   \\
    QBRANDY     & 51\,100   &   2.3   &  63\,100   &   2.1   & 108\,400   &   4.3   &  90\,432   &  24.6   \\
    QCAPRI      &687\,700   &  30.0   &3\,577\,300 & 121.9   &4\,539\,300 & 161.2   &1\,939\,104 & 524.5   \\
    QE226       & 27\,000   &   1.3   &  67\,900   &   2.8   & 115\,400   &   4.2   & 157\,824   &  43.2   \\
    QISRAEL     & 24\,300   &   1.3   &  50\,000   &   2.3   &  65\,200   &   2.9   & 108\,096   &  30.3   \\
    QSC205      & 14\,700   &   1.3    & 14\,800   &    1.2  &  21\,100   &   1.5   &  58\,656   &  16.1   \\
    QSCAGR25    & 20\,100   &   1.1   &  59\,200   &   2.7   &  61\,800   &   2.4   &  62\,400   &  17.5   \\
    QSCFXM3     &277\,300   &  12.7   &2\,123\,400 &  73.4   &1\,284\,100 &  47.0   & 428\,544   & 119.7   \\
    QSEBA       &111\,900   &   6.0   & 160\,500   &   8.4   & 165\,200   &   8.1   & 623\,520   & 166.8   \\
    QSHARE1B    & 72\,200   &   3.5   &2\,020\,400 &  67.3   &2\,324\,300 &  83.1   & 312\,672   &  81.6   \\
    QSHIP12L    & 70\,600   &  4.5    &  63\,200   &  2.6    & 75\,500    &   3.5   &  135\,936  & 38.4    \\
    QSIERRA     &  8\,800   &   0.5   &  40\,400   &   1.3   &  39\,800   &   1.5   & 239\,328   &  64.4   \\
    \midrule
    SGM10        & 46\,562   &   5.0   & 129\,068   &  11.0   &162\,336  &  12.8   & 177\,259   &  55.3   \\
    \bottomrule
  \end{tabular}
\end{table}

\begin{figure}[htp]
  \centering
  \begin{subfigure}[b]{0.48\textwidth}
    \centering
    \includegraphics[width=\textwidth]{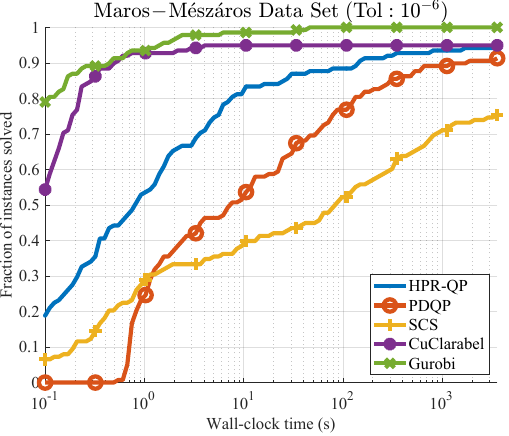}
  \end{subfigure}
  \hfill
  \begin{subfigure}[b]{0.48\textwidth}
    \centering
    \includegraphics[width=\textwidth]{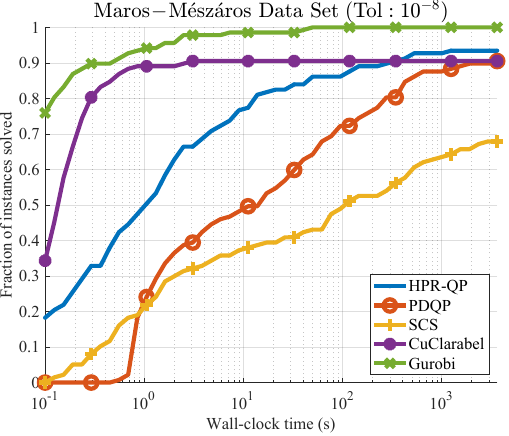}
  \end{subfigure}
  \centering
  \caption{Absolute performance profiles of solvers on the 137 Maros-M\'esz\'aros QP instances}\label{Fig:MM-1e61e8}
\end{figure}

We compare {HPR-QP} with PDQP, SCS, CuClarabel, and Gurobi on the Maros--M\'esz\'aros dataset—a relatively small-scale benchmark—with results summarized in Table~\ref{tab:MM-1e6-1e8} and visualized in Figure~\ref{Fig:MM-1e61e8}. Among all first-order methods tested, {HPR-QP} delivers the best overall performance: it solves more instances than PDQP at both $10^{-6}$ and $10^{-8}$ tolerances and significantly outperforms SCS in both robustness and efficiency. Specifically, {HPR-QP} achieves substantially lower SGM10 in runtime—approximately $3.1\times$ faster at $10^{-6}$ and $3.4\times$ faster at $10^{-8}$—compared to PDQP. While CuClarabel achieves a lower SGM10 in runtime, it solves fewer problems than {HPR-QP} at the $10^{-8}$ tolerance.

\begin{table}[htp]
  \centering
  \caption{Numerical performance on 137 instances of the Maros-M\'esz\'aros data set (Tol.\ $10^{-6}$ and $10^{-8}$).}
  \label{tab:MM-1e6-1e8}
  \renewcommand{\arraystretch}{1.0}
  \begin{tabularx}{\textwidth}{l *{4}{>{\centering\arraybackslash}X}}
    \toprule
    \multirow{2}{*}{Solver}
      & \multicolumn{2}{c}{$10^{-6}$}
      & \multicolumn{2}{c}{$10^{-8}$} \\
    \cmidrule(lr){2-3}\cmidrule(lr){4-5}
      & SGM10 (Time) & Solved & SGM10 (Time) & Solved \\
    \midrule
    HPR-QP       & 10.5  & 129    & 12.6  & 128    \\
    PDQP         & 33.1  & 125    & 42.5  & 124    \\
    SCS          & 126.0    & 103     & 165.0 &  93    \\
    \midrule
    CuClarabel   &  3.7  & 130    &  7.8  & 124    \\
    Gurobi       &  0.4  & 137    &  1.2  & 135    \\
    \bottomrule
  \end{tabularx}
\end{table}

\subsection{Numerical Results on Synthetic CQP Problems}\label{numerical-Sy-QP}

\begin{figure}[H]
  \centering
  \begin{subfigure}[b]{0.48\textwidth}
    \centering
    \includegraphics[width=\textwidth]{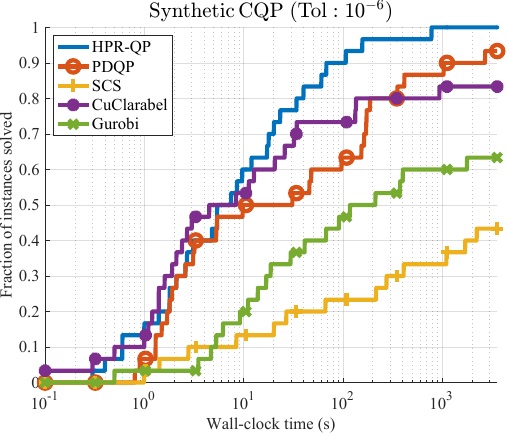}
  \end{subfigure}
  \hfill
  \begin{subfigure}[b]{0.48\textwidth}
    \centering
    \includegraphics[width=\textwidth]{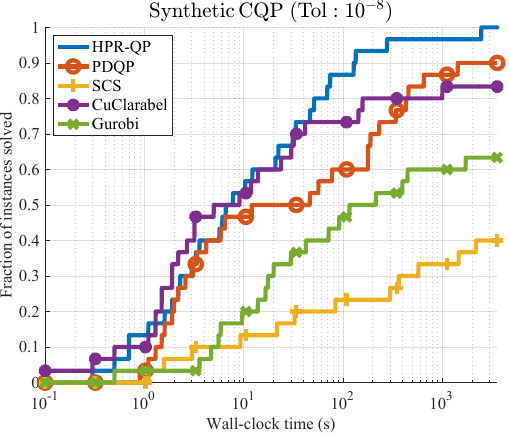}
  \end{subfigure}
  \caption{Absolute performance profiles of tested solvers on 30 synthetic CQP instances.}\label{Fig:Random-1e61e8}
\end{figure}

To complement our evaluation on the relatively small Maros‑Mészáros dataset, we also generated 30 large-scale synthetic CQP instances across six problem classes following Section~8 of \citep{stellato2020osqp}. Table~\ref{tab:instance-sizes} in Appendix~\ref{Append-numerical-result} details their dimensions and sparsity. The numerical performance of all tested solvers is summarized in Table~\ref{tab:random-results-1e6-1e8} and visualized in Figure~\ref{Fig:Random-1e61e8}. At the tighter $10^{-8}$ accuracy level, {HPR‑QP} successfully solves all synthetic instances within one hour. Its SGM10 in runtime is approximately $2.2\times$ faster than the best-performing second-order solver, CuClarabel. When compared with first-order methods, {HPR‑QP} still leads: as shown in Figure~\ref{Fig:Random-1e61e8}, it solves about 90\% of the instances in approximately 100~seconds, whereas the next-best first-order solver, PDQP, requires nearly 1,000~seconds to achieve the same success rate. These results clearly demonstrate the superior efficiency and scalability of {HPR‑QP}.

\begin{table}[H]
  \centering
  \caption{Numerical performance on 30 randomly generated CQP instances (Tol.\ $10^{-6}$ and $10^{-8}$).}
  \label{tab:random-results-1e6-1e8}
  \renewcommand{\arraystretch}{1.2}
  \setlength{\tabcolsep}{10pt}
  \begin{tabular*}{\textwidth}{l @{\extracolsep{\fill}} cc  cc}
    \toprule
    Solver       & \multicolumn{2}{c}{$10^{-6}$} & \multicolumn{2}{c}{$10^{-8}$} \\
    \cmidrule(lr){2-3} \cmidrule(lr){4-5}
                 & SGM10 (Time) & Solved & SGM10 (Time) & Solved \\
    \midrule
    HPR-QP       &  14.3  & 30     &  19.6  & 30     \\
    PDQP         &  51.9  & 28     &  63.8  & 27     \\
    SCS          & 781.5  & 13     & 847.7  & 12     \\
    \midrule
    CuClarabel   &  41.4  & 25     &  43.1  & 25     \\
    Gurobi       & 238.2  & 19     & 242.8  & 19     \\
    \bottomrule
  \end{tabular*}
\end{table}

\subsection{Numerical Results on Lasso Instances}\label{numerical-Lasso}
Consider the following Lasso problem:
\begin{equation}
\label{eq:lasso-org}
\min_{x \in \mathbb{R}^q} \left\{\frac{1}{2}\|\widehat{A} x - \widehat{b}\|^2 + \lambda \|x\|_1 \right\},
\end{equation}
where \( \widehat{A} \in \mathbb{R}^{p \times q} \), \( \widehat{b} \in \mathbb{R}^p \), and \( \lambda > 0 \) is the regularization parameter. This problem can be equivalently reformulated as the following CQP:
\begin{equation}
\label{eq:lasso-convexqp}
\min_{(x,s,t) \in \mathbb{R}^q \times  \mathbb{R}^p \times\mathbb{R}^q} 
\left\{\frac{1}{2} \|s\|^2 + \lambda \sum_{i=1}^q t_i \ \middle| \ s = \widehat{A}x -  \widehat{b},\ -t_i \leq x_i \leq t_i,\ i = 1, \ldots, q \right\}.
\end{equation}
We evaluate performance on 11 Lasso instances from the UCI Machine Learning Repository (as used in \cite{li2018highly}), originally sourced from the LIBSVM datasets \citep{chang2011libsvm}. Table~\ref{tab:uci-info} in Appendix~\ref{Append-numerical-result} reports detailed problem dimensions and sparsity levels of the data matrix \(\widehat{A}\). Numerical results are shown in Table~\ref{tab:uci-1e6-1e8}. In terms of SGM10 in runtime at the \(10^{-8}\) tolerance, {HPR-QP} significantly outperforms both the first-order solver PDQP and the second-order solver CuClarabel. Specifically, {HPR-QP} is approximately \(1.98\times\) faster than CuClarabel and over \(12\times\) faster than PDQP, clearly demonstrating its superiority in both efficiency and scalability across solver categories.

\begin{table}[H]
\centering
\caption{Numerical performance on 11 Lasso instances (Tol.\ $10^{-6}$ and  $10^{-8}$, $\lambda=10^{-3}\|\widehat{A}^{*}\widehat{b}\|_{\infty}$).
{HPR‑QP} solves the original Lasso formulation~\eqref{eq:lasso-org}, while other solvers are applied to its CQP reformulation~\eqref{eq:lasso-convexqp}. Abbreviations: ‘T’ = time‐limit, ‘F’ = failure (e.g. unbounded or infeasible).}
\label{tab:uci-1e6-1e8}
\renewcommand{\arraystretch}{1.1}
\begin{tabularx}{\textwidth}{l *{10}{>{\centering\arraybackslash}X}}
\toprule
 & \multicolumn{2}{c}{HPR-QP}
 & \multicolumn{2}{c}{PDQP}
 & \multicolumn{2}{c}{SCS}
 & \multicolumn{2}{c}{CuClarabel}
 & \multicolumn{2}{c}{Gurobi} \\
\cmidrule(lr){2-3}
\cmidrule(lr){4-5}
\cmidrule(lr){6-7}
\cmidrule(lr){8-9}
\cmidrule(lr){10-11}
Instance 
 & $10^{-6}$ & $10^{-8}$ 
 & $10^{-6}$ & $10^{-8}$ 
 & $10^{-6}$ & $10^{-8}$ 
 & $10^{-6}$ & $10^{-8}$ 
 & $10^{-6}$ & $10^{-8}$ \\
\midrule
abalone7       & 4.2   & 10.5   & 248.6  & 372.5  & T   & T    & 23.2  & 24.4  & 109.1  & 127.3  \\
bodyfat7       & 1.0   & 1.2    & 27.8   & 33.3   & T   & T    & 2.0   & 2.2   & 29.8  & 30.8  \\
E2006.test     & 0.1   & 0.2    & 1.2    & 1.3    & T   & T    & 10.5  & 15.4  & 8.7   & 9.0   \\
E2006.train    & 0.4   & 0.7    & 1.8    & 1.9    & F   & F    & 114.0 & 116.0 & 272.8 & 277.8 \\
housing7       & 11.0  & 22.6   & 83.5   & 123.3  & T   & T    & 5.5   & 5.7   & 123.5 & 125.9 \\
log1p.E2006.test  & 5.0   & 7.0    & 1094.6 & 1416.9 & T   & T    & 183.0 & 196.0 & 107.1 & 137.0 \\
log1p.E2006.train & 13.9  & 17.3   & 2475.7 & 2983.2 & T   & T    & 335.0 & 361.0 & 593.7 & 878.8 \\
mpg7           & 0.4   & 0.6    & 11.6   & 18.1   & 1200.0   & 2000.0    & 0.2   & 0.3   & 1.2   & 1.2   \\
pyrim5         & 35.4  & 49.1   & 298.1  & 410.6  & T   & T    & 3.4   & 3.5   & 35.7  & 35.9  \\
space\_ga9     & 0.5   & 0.6    & 34.9   & 62.7   & 988.0   & 1210.0    & 6.1   & 6.7   & 33.4  & 38.1  \\
triazines4     & 130.7 & 401.3  & 2546.0 & 3533.3 & T   & T    & 25.5  & 26.0  & 455.2 & 843.1 \\
\midrule
SGM10 (Time)          & 8.5   & 13.2   & 124.2  & 161.8  & 2898.0   & 3091.0    & 24.5  & 26.1  & 78.0  & 91.2  \\
\bottomrule
\end{tabularx}
\end{table}

To further evaluate scalability, we generated 12 larger Lasso instances following Appendix~A.5 of \citep{stellato2020osqp}. Numerical results are summarized in Table~\ref{tab:lasso-random-1e6-1e8}. For the relatively small instances (\(p = 10^4\)), {HPR-QP} solves the problems in under 0.5 seconds—faster than the second-best solver, CuClarabel, which takes about 2.5–3.8 seconds. On larger instances (\(p \ge 2 \times 10^5\)), both CuClarabel and Gurobi fail to solve the problems due to time limits or memory exhaustion. In contrast, {HPR-QP} remains robust and efficient. Against the first-order solver {PDQP}, HPR-QP consistently achieves speedups of approximately \(3.9\times\) to \(6.3\times\) in runtime at the \(10^{-8}\) tolerance for these large-scale cases.

\begin{table}[H]
\centering
\caption{Numerical performance on randomly generated Lasso instances  
         (Tol.\ $10^{-6}$ and $10^{-8}$, 
         $\lambda=\frac{1}{5}\|A^*b\|_\infty$).
         {HPR‑QP} solves the original Lasso formulation~\eqref{eq:lasso-org}, while other solvers are applied to its CQP reformulation~\eqref{eq:lasso-convexqp}.
         ‘T’ = time‐limit, ‘M’ = out‐of‐memory.}
\label{tab:lasso-random-1e6-1e8}
\renewcommand{\arraystretch}{1.1}
\begin{tabularx}{\textwidth}{
  >{\centering\arraybackslash}p{1.5cm}
  >{\centering\arraybackslash}p{1.5cm}
  >{\centering\arraybackslash}X >{\centering\arraybackslash}X   
  >{\centering\arraybackslash}X >{\centering\arraybackslash}X   
  >{\centering\arraybackslash}X >{\centering\arraybackslash}X   
  >{\centering\arraybackslash}X >{\centering\arraybackslash}X   
  >{\centering\arraybackslash}X >{\centering\arraybackslash}X   
}
\toprule
\multicolumn{2}{c}{{Size}}
  & \multicolumn{2}{c}{{HPR-QP}}
  & \multicolumn{2}{c}{{PDQP}}
  & \multicolumn{2}{c}{{SCS}}
  & \multicolumn{2}{c}{{CuClarabel}}
  & \multicolumn{2}{c}{{Gurobi}} \\
\cmidrule(lr){1-2}
\cmidrule(lr){3-4}
\cmidrule(lr){5-6}
\cmidrule(lr){7-8}
\cmidrule(lr){9-10}
\cmidrule(lr){11-12}
$p$ & $q$
  & $10^{-6}$ & $10^{-8}$
  & $10^{-6}$ & $10^{-8}$
  & $10^{-6}$ & $10^{-8}$
  & $10^{-6}$ & $10^{-8}$
  & $10^{-6}$ & $10^{-8}$ \\
\midrule
$10^4$         & $5\times10^5$  &  0.5 & 0.5  &  6.6  & 8.7   &  T    &  T    &  2.3  & 2.5   & 22.4 & 23.3  \\
$10^4$         & $10^6$         &  0.2 & 0.3  & 12.2  &17.6   &  T    &  T    &  3.5  & 3.8   & 38.5 & 40.3  \\
$2\times10^5$  & $5\times10^6$  &  5.8 & 9.0  & 38.3  &51.4   &  T    &  T    &   M    &  M    &  T   &  T    \\
$2\times10^5$  & $10^7$         & 16.5 &24.7  &119.1  &155.5  &  T    &  T    &   M    &  M    &  T   &  T    \\
$4\times10^5$  & $5\times10^6$  & 10.7 &14.8  & 60.1  &76.9   &  T    &  T    &   M    &  M    &  M   &  M    \\
$4\times10^5$  & $10^7$         & 32.0 &47.8  &182.8  &238.8  &  T    &  T    &   M    &  M    &  M   &  M    \\
$6\times10^5$  & $5\times10^6$  & 14.2 &21.9  & 76.6  &86.7   &  T    &  T    &   M    &  M    &  M   &  M    \\
$6\times10^5$  & $10^7$         & 47.6 &71.2  &257.5  &328.7  &  T    &  T    &   M    &  M    &  M   &  M    \\
$8\times10^5$  & $5\times10^6$  & 16.7 &22.8  & 93.0  &120.4  &  T    &  T    &   M    &  M    &  M   &  M    \\
$8\times10^5$  & $10^7$         & 63.1 &89.2  &336.2  &420.5  &  T    &  T    &   M    &  M    &  M   &  M    \\
$10^6$         & $5\times10^6$  & 20.7 &25.8  &122.5  &152.3  &  T    &  T    &   M    &  M    &  M   &  M    \\
$10^6$         & $10^7$         & 78.7 &117.7 &   M    &   M    &  T    &  T    &   M    &  M    &  M   &  M    \\
\midrule
\multicolumn{2}{c}{{SGM10 (Time)}}          & 18.5   & 25.5  & 113.2  & 139.3  & 3600.0  & 3600.0   & 1401.3 & 1405.3  & 1691.8 & 1701.0  \\
\bottomrule
\end{tabularx}
\end{table}

\subsection{Numerical Results on QAP Relaxation Problems}\label{numerical-QAP}
Given matrices \(  \widehat{A}, \widehat{B} \in \mathbb{S}^d \) (the space of $d \times d$ symmetric matrices), the QAP is defined as
\[
\min_{X \in \mathbb{R}^{d \times d}} \left\{ \langle \operatorname{vec}(X), ( \widehat{B} \otimes  \widehat{A}) \operatorname{vec}(X) \rangle \;\middle|\; X \in \{0,1\}^{d \times d},\; Xe = e,\; X^* e = e \right\},
\]
where $\otimes$ denotes the Kronecker product, $\operatorname{vec}(X)$ is the vectorization of the matrix $X$, and \( e \in \mathbb{R}^d \) is the vector of all ones. As shown in~\cite{anstreicher2001new}, a good lower bound for the above QAP can be obtained by solving the following CQP:
\[
\min_{\operatorname{vec}(X) \in \mathbb{R}^{d^2}} \left\{ \langle  \operatorname{vec}(X), \widehat{Q} \operatorname{vec}(X) \rangle \;\middle|\; (e^* \otimes I_d) \operatorname{vec}(X) = e,\; (I_d \otimes e^*) \operatorname{vec}(X) = e,\; \operatorname{vec}(X) \geq 0 \right\},
\]
where the matrix \(  \widehat{Q} \in \mathbb{S}^{d^2} \) corresponds to the self-adjoint positive semidefinite linear operator \( \mathcal{ \widehat{Q}} \) defined by
\(
\mathcal{ \widehat{Q}}(X) := ( \widehat{A} X  \widehat{B} - S X - X T)\, \forall X \in \mathbb{R}^{d \times d}.
\)
Equivalently, the matrix form is given by
\(
 \widehat{Q} = ( \widehat{B} \otimes  \widehat{A} - I \otimes S - T \otimes I),
\)
with \( S, T \in \mathbb{S}^d \) computed as follows. Let \(  \widehat{A} = V_A D_A V_A^* \) and \(  \widehat{B} = V_B D_B V_B^* \) be the eigenvalue decompositions of \(  \widehat{A} \) and \(  \widehat{B} \), where
\(
D_A = \operatorname{diag}(\alpha_1, \dots, \alpha_d), \, D_B = \operatorname{diag}(\beta_1, \dots, \beta_d),
\)
with \( \alpha_1 \geq \cdots \geq \alpha_d \) and \( \beta_1 \leq \cdots \leq \beta_d \), and where \( \operatorname{diag}(a_1, \dots, a_d) \) denotes a diagonal matrix with entries \( a_1, \dots, a_d \) on the main diagonal. Let \( (\bar{s}, \bar{t}) \) be the solution to the linear program:
\(
\max \left\{ \langle e, s+t\rangle  \;\middle|\; s_i + t_j \leq \alpha_i \beta_j,\; \forall i,j = 1, \dots, d \right\},
\)
which admits a closed-form solution as described in~\cite{anstreicher2001new}. Then,
\(
S = V_A \operatorname{diag}(\bar{s}) V_A^*, T = V_B \operatorname{diag}(\bar{t}) V_B^*.
\)

\begin{figure}[H]
  \centering
  \begin{subfigure}[b]{0.48\textwidth}
    \centering
    \includegraphics[width=\textwidth]{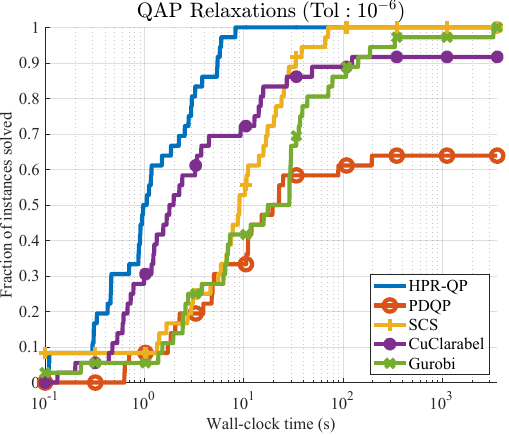}
  \end{subfigure}
  \hfill
  \begin{subfigure}[b]{0.48\textwidth}
    \centering
    \includegraphics[width=\textwidth]{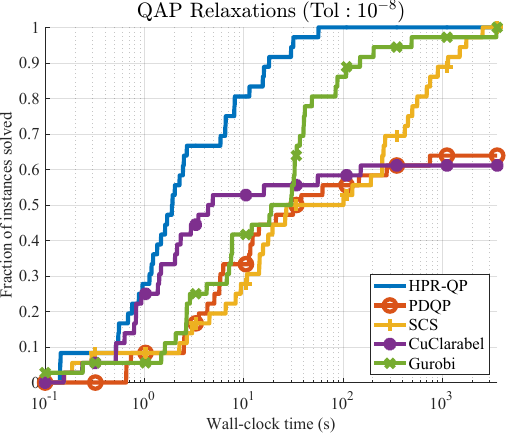}
  \end{subfigure}
  \caption{Absolute performance profiles of tested solvers on 36 instances of QAP relaxations.}\label{Fig:QAP-1e61e8}
\end{figure}

We tested 36 instances with $50\leq d \leq 256$ from the QAPLIB benchmark set~\citep{burkard1997qaplib}.\footnote{The instance {sko90} is excluded, as Gurobi reports it to be non-convex.} The numerical results are summarized in Figure~\ref{Fig:QAP-1e61e8} and Table~\ref{tab:real-QAP-1e6-1e8}. Among the first-order solvers, Table~\ref{tab:real-QAP-1e6-1e8} shows that {HPR-QP} significantly outperforms SCS, the second-best solver in this category, achieving a speedup of approximately \(18.3\times\) at \(10^{-8}\) tolerance and \(6.3\times\) at \(10^{-6}\) tolerance. Figure~\ref{Fig:QAP-1e61e8} further demonstrates that for \(10^{-8}\) accuracy, {HPR-QP} solves around 90\% of the instances within 20 seconds, while SCS requires nearly 1000 seconds to reach the same success rate.

\begin{table}[H]
  \centering
  \caption{Numerical performance on 36 instances of QAP relaxations (Tol.\ $10^{-6}$ and $10^{-8}$).}
  \label{tab:real-QAP-1e6-1e8}
  \renewcommand{\arraystretch}{1.0}
  \begin{tabularx}{\linewidth}{l *{4}{>{\centering\arraybackslash}X}}
    \toprule
    \multirow{2}{*}{Solver}
      & \multicolumn{2}{c}{$10^{-6}$}
      & \multicolumn{2}{c}{$10^{-8}$} \\
    \cmidrule(lr){2-3} \cmidrule(lr){4-5}
      & SGM10 (Time) & Solved & SGM10 (Time) & Solved \\
    \midrule
    HPR‐QP  &  1.8   & 36     &  4.7   & 36     \\
    PDQP               & 124.1  & 23     & 149.4  & 23     \\
    SCS                & 11.3   & 36     & 86.0   & 36     \\
    \midrule
    CuClarabel           & 13.6   & 33     &114.9   & 22     \\
    Gurobi             & 24.8   & 36     & 26.8   & 36     \\
    \bottomrule
  \end{tabularx}
\end{table}

To further evaluate the scalability of the solvers, we generate extremely large-scale synthetic QAP instances. Specifically, we uniformly sample \( d \) points in the unit square \([0,1)^2\), and define the matrix \( \widehat{A}\) as their pairwise Euclidean distance matrix, i.e., \( \widehat{A}_{ij} = \|(x_i, y_i) - (x_j, y_j)\| \). The matrix \( \widehat{B} \) is constructed as a symmetric random matrix with entries drawn uniformly from \([0,1)\) and zero diagonal. Due to the prohibitive size of \( d \) and the limited memory, forming an explicit matrix representation of \( Q \) becomes infeasible. Consequently, we only report the results for {HPR-QP} in Table~\ref{tab:QAP-large}, highlighting its crucial advantage of operating without an explicit form of \( Q \). Notably, {HPR‑QP} successfully solves QAP relaxation instances up to \(d = 8192\), demonstrating its exceptional scalability and efficiency in matrix‑free settings.

\begin{table}[H]
\centering
\caption{The runtime performance of HPR-QP on extremely large-scale QAP relaxations.}
\label{tab:QAP-large}
\renewcommand{\arraystretch}{1.2}
\begin{tabularx}{\linewidth}{l *{3}{>{\centering\arraybackslash}X}}
\toprule
$d$               & $10^{-4}$  & $10^{-6}$   & $10^{-8}$ \\
\midrule
256       & 0.8   & 3.0     & 6.0     \\
512       & 1.3   & 4.0     & 9.3     \\
1024      & 3.6   & 26.3    & 73.3    \\
2048      & 73.0   & 291.0  & 700.0   \\
4096      & 32.0  & 3640.0  & 9580.7  \\
8192      & 188.6 & 61905.6 & 126547.9\\
\bottomrule
\end{tabularx}
\end{table}

\section{Conclusion}\label{sec:5}
In this paper, we proposed {HPR-QP}, a dual HPR method designed for solving large-scale CCQP problems. HPR-QP incorporates adaptive restart and penalty parameter update strategies to enhance convergence and robustness. Extensive numerical experiments on benchmark CCQP data sets demonstrate that {HPR-QP} delivers competitive performance across a wide range of accuracy requirements, particularly excelling in large-scale scenarios where second-order methods face limitations in scalability. Nonetheless, for small-scale problems, second-order methods may still offer superior efficiency due to their fast local convergence. As a potential direction for future work, integrating direct solvers into the {HPR-QP} framework to more efficiently handle the linear systems in subproblems could further improve its performance on smaller instances. Additionally, it would be interesting to explore hybrid strategies that combine {HPR-QP} with second-order methods on GPU, enabling more adaptive solutions across different problem scales.

\bibliographystyle{abbrv}
\bibliography{HPR-QP-ref}

\begin{thebibliography}{10}

\bibitem{anstreicher2001new}
K.~M. Anstreicher and N.~W. Brixius.
\newblock A new bound for the quadratic assignment problem based on convex
  quadratic programming.
\newblock {\em Math. Program.}, 89:341--357, 2001.

\bibitem{applegate2021practical}
D.~Applegate, M.~D{\'\i}az, O.~Hinder, H.~Lu, M.~Lubin, B.~O'Donoghue, and
  W.~Schudy.
\newblock Practical large-scale linear programming using primal-dual hybrid
  gradient.
\newblock In {\em Advances in Neural Information Processing System}, volume~34,
  pages 20243--20257, 2021.

\bibitem{applegate2023faster}
D.~Applegate, O.~Hinder, H.~Lu, and M.~Lubin.
\newblock Faster first-order primal-dual methods for linear programming using
  restarts and sharpness.
\newblock {\em Math. Program.}, 201(1):133--184, 2023.

\bibitem{bezanson2017julia}
J.~Bezanson, A.~Edelman, S.~Karpinski, and V.~B. Shah.
\newblock Julia: A fresh approach to numerical computing.
\newblock {\em SIAM Rev.}, 59(1):65--98, 2017.

\bibitem{bredies2022degenerate}
K.~Bredies, E.~Chenchene, D.~A. Lorenz, and E.~Naldi.
\newblock Degenerate preconditioned proximal point algorithms.
\newblock {\em SIAM J. Optim.}, 32(3):2376--2401, 2022.

\bibitem{burkard1997qaplib}
R.~E. Burkard, S.~E. Karisch, and F.~Rendl.
\newblock {Q}{A}{P}{L}{I}{B}--a quadratic assignment problem library.
\newblock {\em J. Glob. Optim.}, 10:391--403, 1997.

\bibitem{chang2011libsvm}
C.-C. Chang and C.-J. Lin.
\newblock {LIBSVM}: {A} library for support vector machines.
\newblock {\em ACM Trans. Intell. Syst. Technol.}, 2(3):1--27, 2011.

\bibitem{chen2024hpr}
K.~Chen, D.~F. Sun, Y.~Yuan, G.~Zhang, and X.~Zhao.
\newblock {HPR-LP}: {A}n implementation of an {HPR} method for solving linear
  programming.
\newblock {\em arXiv preprint arXiv:2408.12179}, 2024.

\bibitem{chenOptimal2014}
Y.~Chen, G.~Lan, and Y.~Ouyang.
\newblock Optimal primal-dual methods for a class of saddle point problems.
\newblock {\em SIAM J. Optim.}, 24(4):1779--1814, 2014.

\bibitem{chen2024cuclarabel}
Y.~Chen, D.~Tse, P.~Nobel, P.~Goulart, and S.~Boyd.
\newblock {CuClarabel}: {GPU} acceleration for a conic optimization solver.
\newblock {\em arXiv preprint arXiv:2412.19027}, 2024.

\bibitem{deng2024enhanced}
Q.~Deng, Q.~Feng, W.~Gao, D.~Ge, B.~Jiang, Y.~Jiang, J.~Liu, T.~Liu, C.~Xue,
  Y.~Ye, et~al.
\newblock An enhanced alternating direction method of multipliers-based
  interior point method for linear and conic optimization.
\newblock {\em INFORMS J. Comput}, 2024.

\bibitem{DornQPDual60}
W.~S. Dorn.
\newblock Duality in quadratic programming.
\newblock {\em Quart. Appl. Math.}, 18:155--162, 1960/61.

\bibitem{eckstein1992douglas}
J.~Eckstein and D.~P. Bertsekas.
\newblock On the {D}ouglas—{R}achford splitting method and the proximal point
  algorithm for maximal monotone operators.
\newblock {\em Math. Program.}, 55(1):293--318, 1992.

\bibitem{gerald2004applied}
C.~F. Gerald.
\newblock {\em Applied Numerical Analysis}.
\newblock Pearson Education India, 2004.

\bibitem{golub2013matrix}
G.~H. Golub and C.~F. Van~Loan.
\newblock {\em Matrix Computations}.
\newblock JHU press, 2013.

\bibitem{goulart2024clarabel}
P.~J. Goulart and Y.~Chen.
\newblock Clarabel: {An} interior-point solver for conic programs with
  quadratic objectives.
\newblock {\em arXiv preprint arXiv:2405.12762}, 2024.

\bibitem{gurobi}
{Gurobi Optimization, LLC}.
\newblock {Gurobi {O}ptimizer {R}eference {M}anual}, 2024.

\bibitem{halpern1967fixed}
B.~Halpern.
\newblock Fixed points of nonexpanding maps.
\newblock {\em Bull. Am. Math. Soc.}, 73(6):957--961, 1967.

\bibitem{han2018linear}
D.~Han, D.~F. Sun, and L.~Zhang.
\newblock Linear rate convergence of the alternating direction method of
  multipliers for convex composite programming.
\newblock {\em Math. Oper. Res.}, 43(2):622--637, 2018.

\bibitem{huang2024restarted}
Y.~Huang, W.~Zhang, H.~Li, D.~Ge, H.~Liu, and Y.~Ye.
\newblock Restarted primal-dual hybrid conjugate gradient method for
  large-scale quadratic programming.
\newblock {\em arXiv preprint arXiv:2405.16160}, 2024.

\bibitem{cplex2009v12}
IBM.
\newblock {IBM ILOG} {CPLEX} {O}ptimization {S}tudio {CPLEX} {U}ser’s
  {M}anual, 1987.

\bibitem{li2015convergent}
M.~Li, D.~F. Sun, and K.-C. Toh.
\newblock A convergent 3-block semi-proximal {ADMM} for convex minimization
  problems with one strongly convex block.
\newblock {\em Asia-Pac. J. Oper. Res.}, 32(04):1550024, 2015.

\bibitem{li2016schur}
X.~Li, D.~F. Sun, and K.-C. Toh.
\newblock A {S}chur complement based semi-proximal {A}{D}{M}{M} for convex
  quadratic conic programming and extensions.
\newblock {\em Math. Program.}, 155(1-2):333--373, 2016.

\bibitem{li2018highly}
X.~Li, D.~F. Sun, and K.-C. Toh.
\newblock A highly efficient semismooth {Newton} augmented {Lagrangian} method
  for solving {L}asso problems.
\newblock {\em SIAM J. Optim.}, 28(1):433--458, 2018.

\bibitem{li2018qsdpnal}
X.~Li, D.~F. Sun, and K.-C. Toh.
\newblock {Q}{S}{D}{P}{N}{A}{L}: {A} two-phase augmented {L}agrangian method
  for convex quadratic semidefinite programming.
\newblock {\em Math. Program. Comput.}, 10:703--743, 2018.

\bibitem{li2019block}
X.~Li, D.~F. Sun, and K.-C. Toh.
\newblock A block symmetric {G}auss--{S}eidel decomposition theorem for convex
  composite quadratic programming and its applications.
\newblock {\em Math. Program.}, 175:395--418, 2019.

\bibitem{liang2022qppal}
L.~Liang, X.~Li, D.~F. Sun, and K.-C. Toh.
\newblock {Q}{P}{P}{A}{L}: {A} two-phase proximal augmented {L}agrangian method
  for high-dimensional convex quadratic programming problems.
\newblock {\em ACM Trans. Math. Softw.}, 48(3):1--27, 2022.

\bibitem{lieder2021convergence}
F.~Lieder.
\newblock On the convergence rate of the {H}alpern-iteration.
\newblock {\em Optim. Lett.}, 15(2):405--418, 2021.

\bibitem{lin2021admm}
T.~Lin, S.~Ma, Y.~Ye, and S.~Zhang.
\newblock An {A}{D}{M}{M}-based interior-point method for large-scale linear
  programming.
\newblock {\em Optim. Methods Softw.}, 36(2-3):389--424, 2021.

\bibitem{lin2025pdcs}
Z.~Lin, Z.~Xiong, D.~Ge, and Y.~Ye.
\newblock {PDCS}: {A} primal-dual large-scale conic programming solver with
  {GPU} enhancements.
\newblock {\em arXiv preprint arXiv:2505.00311}, 2025.

\bibitem{lions1979splitting}
P.-L. Lions and B.~Mercier.
\newblock Splitting algorithms for the sum of two nonlinear operators.
\newblock {\em SIAM J. Numer. Anal.}, 16(6):964--979, 1979.

\bibitem{lu2023cupdlp}
H.~Lu and J.~Yang.
\newblock cu{P}{D}{L}{P}.jl: {A} {G}{P}{U} implementation of restarted
  primal-dual hybrid gradient for linear programming in {J}ulia.
\newblock {\em arXiv preprint arXiv:2311.12180}, 2023.

\bibitem{lu2023practical}
H.~Lu and J.~Yang.
\newblock A practical and optimal first-order method for large-scale convex
  quadratic programming.
\newblock {\em arXiv preprint arXiv:2311.07710}, 2023.

\bibitem{lu2023cupdlp_c}
H.~Lu, J.~Yang, H.~Hu, Q.~Huangfu, J.~Liu, T.~Liu, Y.~Ye, C.~Zhang, and D.~Ge.
\newblock cu{P}{D}{L}{P}-{C}: {A} strengthened implementation of cu{P}{D}{L}{P}
  for linear programming by {C} language.
\newblock {\em arXiv preprint arXiv:2312.14832}, 2023.

\bibitem{maros1999repository}
I.~Maros and C.~M{\'e}sz{\'a}ros.
\newblock A repository of convex quadratic programming problems.
\newblock {\em Optim. Methods Softw.}, 11(1-4):671--681, 1999.

\bibitem{o2021operator}
B.~O'Donoghue.
\newblock Operator splitting for a homogeneous embedding of the linear
  complementarity problem.
\newblock {\em SIAM J. Optim.}, 31(3):1999--2023, 2021.

\bibitem{o2016conic}
B.~O’donoghue, E.~Chu, N.~Parikh, and S.~Boyd.
\newblock Conic optimization via operator splitting and homogeneous self-dual
  embedding.
\newblock {\em J. Optim. Theory Appl.}, 169:1042--1068, 2016.

\bibitem{pock2011diagonal}
T.~Pock and A.~Chambolle.
\newblock Diagonal preconditioning for first order primal-dual algorithms in
  convex optimization.
\newblock In {\em 2011 International Conference on Computer Vision}, pages
  1762--1769. IEEE, 2011.

\bibitem{rockafellar1970convex}
R.~T. Rockafellar.
\newblock {\em Convex {A}nalysis}.
\newblock Princeton {U}niversity {P}ress, 1970.

\bibitem{ruiz2001scaling}
D.~Ruiz.
\newblock A scaling algorithm to equilibrate both rows and columns norms in
  matrices.
\newblock Technical report, Rutherford Appleton Laboratory, 2001.

\bibitem{sabach2017first}
S.~Sabach and S.~Shtern.
\newblock A first order method for solving convex bilevel optimization
  problems.
\newblock {\em SIAM J. Optim.}, 27(2):640--660, 2017.

\bibitem{stellato2020osqp}
B.~Stellato, G.~Banjac, P.~Goulart, A.~Bemporad, and S.~Boyd.
\newblock {O}{S}{Q}{P}: {A}n operator splitting solver for quadratic programs.
\newblock {\em Math. Program. Comput.}, 12(4):637--672, 2020.

\bibitem{sun2025accelerating}
D.~F. Sun, Y.~Yuan, G.~Zhang, and X.~Zhao.
\newblock Accelerating preconditioned {ADMM} via degenerate proximal point
  mappings.
\newblock {\em SIAM J. Optim.}, 35(2):1165--1193, 2025.

\bibitem{WolfeDual1961}
P.~Wolfe.
\newblock A duality theorem for non-linear programming.
\newblock {\em Quart. Appl. Math.}, 19:239--244, 1961.

\bibitem{yang2025accelerated}
B.~Yang, X.~Zhao, X.~Li, and D.~F. Sun.
\newblock An accelerated proximal alternating direction method of multipliers
  for optimal decentralized control of uncertain systems.
\newblock {\em J. Optim. Theory Appl.}, 204(1):9, 2025.

\bibitem{zhang2025hot}
G.~Zhang, Z.~Gu, Y.~Yuan, and D.~F. Sun.
\newblock H{O}{T}: {A}n efficient {H}alpern accelerating algorithm for optimal
  transport problems.
\newblock {\em IEEE Trans. Pattern Anal. Mach. Intell. XXX, in print.}, 2025.

\bibitem{zhang2022efficient}
G.~Zhang, Y.~Yuan, and D.~F. Sun.
\newblock An efficient {H}{P}{R} algorithm for the {W}asserstein barycenter
  problem with ${O}$({D}im({P})/$\varepsilon$) computational complexity.
\newblock {\em arXiv preprint arXiv:2211.14881}, 2022.

\end{thebibliography}

\newpage
\appendix

\section{An HPR Method for the Primal Form of CCQP}\label{appendix-primal}
\subsection{An HPR Method Applied to the Primal Reformulation \texorpdfstring{\eqref{QP-primal-s}}{(1.3)}%
}\label{appendix-primal-1}
Recall the following primal reformulation of CCQP problems:
\begin{equation}\label{QP-primal-4}
   \min _{(x,s) \in   \mathbb{R}^n \times \mathbb{R}^{m}}\left\{\frac{1}{2}\langle x, Q x\rangle +\langle c,x\rangle +\phi(x)+\delta_\cK(s)    \mid Ax=s \right\}.
   \end{equation}
 Given $\sigma>0$, we define the augmented Lagrangian function $L_\sigma(x,s;y)$ associated with problem \eqref{QP-primal-4} for any $(x,s,y) \in \mathbb{R}^{n} \times   \mathbb{R}^{m} \times \mathbb{R}^m$ as follows:
$$
L_\sigma(x,s;y)=\frac{1}{2}\langle x, Qx\rangle+\langle c,x\rangle+\phi(x)+\delta_{\cK}(s) + \langle s-Ax,  y \rangle  +\frac{\sigma}{2}\left\|Ax-s\right\|^2.
$$
For the sake of notational simplicity, we represent the tuple $(x,s,y)$ as $u$, and we denote the set $ \mathbb{R}^{n} \times  \mathbb{R}^{m} \times \mathbb{R}^m$ by $\mathcal{U}$. Then, an HPR method \citep{sun2025accelerating} for solving problem \eqref{QP-primal-4} is presented in Algorithm \ref{alg:acc-pADMM-QP-p-1}: 
\begin{algorithm}[H] 
	\caption{An HPR method for solving the primal reformulation \eqref{QP-primal-4}}
	\label{alg:acc-pADMM-QP-p-1}
	\begin{algorithmic}[1] 
	\State {Input:  Choose $u^{0}=(x^0,s^0,y^0)\in    \mathbb{R}^{n}  \times \mathbb{R}^{m} \times   \mathbb{R}^{m}$. Set parameters $\sigma > 0$.  Let $\cS_x=\lambda_QI_n-Q+\sigma(\lambda_AI_n-A^*A)$. Denote $u=(x,s,y)$ and $\bu=(\bx,\bs,\by)$.}
    \For{\( k = 0, 1, 2, \ldots \)}
    \State {Step 1. $\displaystyle \bs^{k+1}=\underset{s \in \mathbb{R}^{m}}{\arg \min }\left\{L_\sigma(x^k,s;y^k)\right\}$;} 
    \State {Step 2. $\displaystyle  \by^{k+1}={y}^k+\sigma (\bs^{k+1}-Ax^k) $;}
    \State {Step 3. $\displaystyle \bx^{k+1}=\underset{x \in \mathbb{R}^{n}}{\arg \min }\left\{L_\sigma(x,\bs^{k+1};\by^{k+1})+\frac{1}{2}\|x-x^k\|^2_{\cS_x}\right\}$;}
    \State {Step 4. $\displaystyle \hat{u}^{k+1}= 2 \bu^{k+1}-u^k$;}
	\State {Step 5. $\displaystyle u^{k+1}=\frac{1}{k+2}u^{0}+\frac{k+1}{k+2}\hat{u}^{k+1}$;}
    \EndFor
	\end{algorithmic}
\end{algorithm}

\subsection{An HPR Method Applied to the Primal Reformulation
  \texorpdfstring{\eqref{QP-primal-sv}}{(1.4)}%
}\label{appendix-primal-2}

Recall the following primal reformulation of CCQP problems:
\begin{equation}\label{QP-primal-5}
   \min _{(x,s,v) \in   \mathbb{R}^n\times \mathbb{R}^{m}  \times\mathbb{R}^{n} }\left\{\frac{1}{2}\langle v, Q v\rangle +\langle c,x\rangle +\phi(x)+\delta_\cK(s)    \mid Ax=s, x=v \right\}.
\end{equation} Given $\sigma>0$, we define the augmented Lagrangian function $L_\sigma(x,s,v;y)$ associated with problem \eqref{QP-primal-5} for any $(x,s,v,y,t) \in \mathbb{R}^{n} \times   \mathbb{R}^{m}\times   \mathbb{R}^{n} \times \mathbb{R}^m \times \mathbb{R}^n$ as follows:
$$
L_\sigma(x,s,v;y,t)=\frac{1}{2}\langle v, Qv\rangle+\langle c,x\rangle+\phi(x)+\delta_{\cK}(s) + \langle s-Ax,  y \rangle   + \langle v-x, t \rangle +\frac{\sigma}{2}\left\|Ax-s\right\|^2+\frac{\sigma}{2}\|x-v\|^2.
$$
For the sake of notational simplicity, we represent the tuple $(x,s,v,y,t)$ as $u$, and we denote the set $ \mathbb{R}^{n} \times  \mathbb{R}^{m} \times \mathbb{R}^n\times \mathbb{R}^m\times \mathbb{R}^n$ by $\mathcal{U}$. Then, an HPR method \citep{sun2025accelerating} for solving problem \eqref{QP-primal-5} can be presented in Algorithm \ref{alg:acc-pADMM-QP-p-2}: 
\begin{algorithm}[H] 
	\caption{An HPR method for solving the primal reformulation \eqref{QP-primal-5}}
	\label{alg:acc-pADMM-QP-p-2}
	\begin{algorithmic}[1] 
	\State {Input:  Choose $u^{0}=(x^0,s^0,v^{0},y^0,t^{0})\in \mathcal{U}$. Set parameters $\sigma > 0$. Let $\cS_v= \lambda_QI_n-Q$ and $\cS_x= \sigma(\lambda_AI_n-A^*A)$. Denote $u=(x,s,v,y,t)$ and $\bu=(\bx,\bs,\bv,\by,\bar t)$.  }
    \For{\( k = 0, 1, 2, \ldots \)} 
    \State {Step 1. $ (\bs^{k+1},\bv^{k+1})=\displaystyle \underset{(s,v) \in \mathbb{R}^{m}\times \mathbb{R}^{n}}{\arg \min }\left\{L_\sigma(x^k,s,v;y^k,t^k)+\frac{1}{2}\|v-v^k\|_{\cS_v}^2\right\} $;} 
    \State {Step 2.1. $\displaystyle  \by^{k+1}={y}^k+\sigma (\bs^{k+1}-Ax^k) $;}
    \State {Step 2.2.  $\displaystyle  \bar t^{k+1}={t}^k+\sigma (\bv^{k+1}-x^k) $;}   
    \State {Step 3. $\displaystyle \bx^{k+1}=\underset{x \in \mathbb{R}^{n}}{\arg \min }\left\{L_\sigma(x,\bs^{k+1},\bv^{k+1};\by^{k+1},\bar t^{k+1})+\frac{1}{2}\|x-x^k\|^2_{\cS_x}\right\}$;}
    \State {Step 4. $\displaystyle \hat{u}^{k+1}= 2 \bu^{k+1}-u^k$;}
	\State {Step 5. $\displaystyle u^{k+1}=\frac{1}{k+2}u^{0}+\frac{k+1}{k+2}\hat{u}^{k+1}$;}
    \EndFor
	\end{algorithmic}
\end{algorithm}

\section{Additional Numerical Results}\label{Append-numerical-result}

\begin{figure}[htp]
  \centering
  \begin{subfigure}[b]{0.48\textwidth}
    \centering
    \renewcommand{\arraystretch}{1.95}
    \begin{tabularx}{\linewidth}{l *{2}{>{\centering\arraybackslash}X}}
      \toprule
      Solver     & SGM10 (Time) & Solved \\
      \midrule
      HPR-QP     &  7.4  & 130    \\
      PDQP       & 11.3  & 133    \\
      SCS        & 38.1  & 123    \\
      \midrule
      CuClarabel &  0.7  & 136    \\
      Gurobi     &  0.4  & 137    \\
      \bottomrule
    \end{tabularx}
    \subcaption{SGM10 and number solved}%
    \label{tab:MM-1e4}
  \end{subfigure}%
  \hfill
  \begin{subfigure}[b]{0.48\textwidth}
    \centering
    \includegraphics[width=\linewidth]{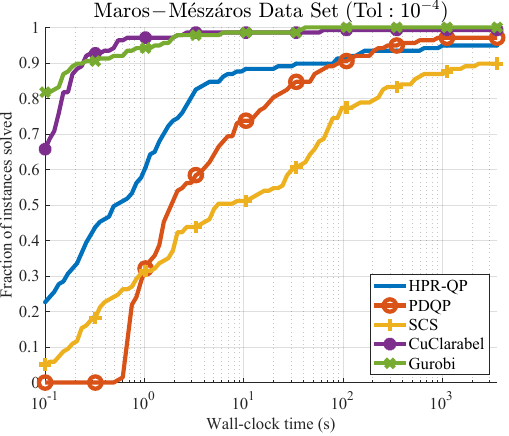}
    \subcaption{Absolute performance profiles}%
    \label{fig:MM-1e4}
  \end{subfigure}

  \caption{Numerical performance of tested solvers on 137 instances of the Maros-M\'esz\'aros data set (Tol.\ $10^{-4}$).}
  \label{fig:MM-1e4-combined}
\end{figure}

\begin{table}[htp]
\centering
\small
\caption{Problem dimensions and sparsity of matrices \(A\) and \(Q\) in the synthetic CQP instances.}
\label{tab:instance-sizes}
\renewcommand{\arraystretch}{1.1}
\begin{tabularx}{\textwidth}{l *{4}{>{\centering\arraybackslash}X}}
\toprule
Instance            & Rows                  & Cols                  & nnz($A$)               & nnz($Q$)              \\
\midrule
random\_1         & $5.00\times10^{5}$    & $5.00\times10^{4}$    & $1.00\times10^{6}$     & $2.50\times10^{5}$    \\
random\_2        & $1.00\times10^{6}$    & $1.00\times10^{5}$    & $2.00\times10^{6}$     & $5.01\times10^{5}$    \\
random\_3        & $5.00\times10^{6}$    & $5.00\times10^{5}$    & $1.00\times10^{7}$     & $2.50\times10^{6}$    \\
random\_4       & $1.00\times10^{7}$    & $1.00\times10^{6}$    & $2.00\times10^{7}$     & $5.00\times10^{6}$    \\
random\_5       & $5.00\times10^{7}$    & $5.00\times10^{6}$    & $1.00\times10^{8}$     & $2.50\times10^{7}$    \\
equality\_1        & $5.00\times10^{3}$    & $5.00\times10^{3}$    & $1.99\times10^{5}$     & $6.80\times10^{6}$    \\
equality\_2       & $5.00\times10^{3}$    & $1.00\times10^{4}$    & $2.00\times10^{5}$     & $1.47\times10^{7}$    \\
equality\_3       & $1.00\times10^{4}$    & $2.00\times10^{4}$    & $4.00\times10^{5}$     & $3.07\times10^{7}$    \\
equality\_4       & $2.50\times10^{4}$    & $5.00\times10^{4}$    & $1.00\times10^{6}$     & $7.87\times10^{7}$    \\
equality\_5      & $5.00\times10^{4}$    & $1.00\times10^{5}$    & $2.00\times10^{6}$     & $1.59\times10^{8}$    \\
control\_1        & $2.20\times10^{3}$    & $3.20\times10^{3}$    & $6.02\times10^{5}$     & $4.24\times10^{4}$    \\
control\_2        & $5.50\times10^{3}$    & $8.00\times10^{3}$    & $3.76\times10^{6}$     & $2.56\times10^{5}$    \\
control\_3         & $1.10\times10^{4}$    & $1.60\times10^{4}$    & $1.50\times10^{7}$     & $1.01\times10^{6}$    \\
control\_4       & $1.65\times10^{4}$    & $2.40\times10^{4}$    & $3.38\times10^{7}$     & $2.27\times10^{6}$    \\
control\_5         & $2.20\times10^{4}$    & $3.20\times10^{4}$    & $6.00\times10^{7}$     & $4.02\times10^{6}$    \\
portfolio\_1      & $4.01\times10^{2}$    & $4.04\times10^{4}$    & $8.04\times10^{6}$     & $4.04\times10^{4}$    \\
portfolio\_2      & $5.01\times10^{2}$    & $5.05\times10^{4}$    & $1.26\times10^{7}$     & $5.05\times10^{4}$    \\
portfolio\_3      & $6.01\times10^{2}$    & $6.06\times10^{4}$    & $1.81\times10^{7}$     & $6.06\times10^{4}$    \\
portfolio\_4      & $7.01\times10^{2}$    & $7.07\times10^{4}$    & $2.46\times10^{7}$     & $7.07\times10^{4}$    \\
portfolio\_5      & $8.01\times10^{2}$    & $8.08\times10^{4}$    & $3.21\times10^{7}$     & $8.08\times10^{4}$    \\
huber\_1           & $5.00\times10^{5}$    & $1.51\times10^{6}$    & $1.60\times10^{6}$     & $5.00\times10^{5}$    \\
huber\_2          & $1.00\times10^{6}$    & $3.01\times10^{6}$    & $3.20\times10^{6}$     & $1.00\times10^{6}$    \\
huber\_3          & $2.00\times10^{6}$    & $6.02\times10^{6}$    & $6.40\times10^{6}$     & $2.00\times10^{6}$    \\
huber\_4          & $5.00\times10^{6}$    & $1.51\times10^{7}$    & $1.60\times10^{7}$     & $5.00\times10^{6}$    \\
huber\_5         & $1.00\times10^{7}$    & $3.01\times10^{7}$    & $3.20\times10^{7}$     & $1.00\times10^{7}$    \\
svm\_1            & $1.00\times10^{6}$    & $1.01\times10^{6}$    & $1.40\times10^{6}$     & $1.00\times10^{4}$    \\
svm\_2            & $2.00\times10^{6}$    & $2.02\times10^{6}$    & $2.80\times10^{6}$     & $2.00\times10^{4}$    \\
svm\_3            & $5.00\times10^{6}$    & $5.05\times10^{6}$    & $7.00\times10^{6}$     & $5.00\times10^{4}$    \\
svm\_4           & $1.00\times10^{7}$    & $1.01\times10^{7}$    & $1.40\times10^{7}$     & $1.00\times10^{5}$    \\
svm\_5           & $2.00\times10^{7}$    & $2.02\times10^{7}$    & $2.80\times10^{7}$     & $2.00\times10^{5}$    \\
\bottomrule
\end{tabularx}
\end{table}

\begin{figure}[htp]
  \centering
  \begin{subfigure}[b]{0.48\textwidth}
    \centering
    \renewcommand{\arraystretch}{1.95}
    \setlength{\tabcolsep}{10pt}
    \begin{tabularx}{\linewidth}{l *{2}{>{\centering\arraybackslash}X}}
      \toprule
      Solver       & SGM10 (Time)  & Solved \\
      \midrule
      HPR‐QP       &   8.6 & 30     \\
      PDQP         &  42.6 & 28     \\
      SCS          & 648.9 & 13     \\
      CuClarabel   &  39.1 & 25     \\
      Gurobi       & 230.4 & 20     \\
      \bottomrule
    \end{tabularx}
    \subcaption{SGM10 and number of problems solved}%
    \label{tab:random-1e4}
  \end{subfigure}%
  \hfill
  \begin{subfigure}[b]{0.48\textwidth}
    \centering
    \includegraphics[width=\linewidth]{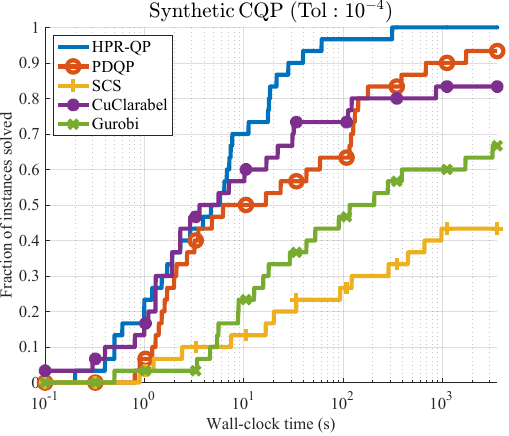}
    \subcaption{Absolute performance profiles}%
    \label{fig:random-1e4}
  \end{subfigure}

  \caption{Numerical performance of tested solvers on 30 synthetic CQP problems (Tol.\ $10^{-4}$).}
  \label{fig:random-1e4-combined}
\end{figure}

\begin{table}[htp]
\centering
\caption{Dimensions and number of nonzeros in $\widehat{A}$ for UCI Lasso instances.}
\label{tab:uci-info}
\renewcommand{\arraystretch}{1.2}
\begin{tabular}{l r r r}
\toprule
Instance             & \multicolumn{1}{c}{$p$}    & \multicolumn{1}{c}{$q$}       & \multicolumn{1}{c}{$\mathrm{nnz}(\widehat A)$} \\
\midrule
abalone7             & 4\,177                     & 6\,435                        & 22\,873\,422                                    \\
bodyfat7             & 252                        & 116\,280                      & 29\,302\,560                                    \\
E2006.test           & 3\,308                     & 150\,358                      & 4\,559\,533                                     \\
E2006.train          & 16\,087                    & 150\,360                      & 19\,971\,015                                    \\
housing7             & 506                        & 77\,520                       & 39\,225\,120                                    \\
log1p.E2006.test     & 3\,308                     & 4\,272\,226                   & 22\,474\,250                                    \\
log1p.E2006.train    & 16\,087                    & 4\,272\,227                   & 96\,731\,839                                    \\
mpg7                 & 392                        & 3\,432                        & 1\,174\,932                                     \\
pyrim5               & 74                         & 201\,376                      & 8\,054\,057                                     \\
space\_ga9           & 3\,107                     & 5\,005                        & 15\,550\,535                                    \\
triazines4           & 186                        & 635\,376                      & 77\,638\,169                                    \\
\bottomrule
\end{tabular}
\end{table}

\begin{table}[htp]
\centering
\caption{Numerical performance on 11 Lasso instances (Tol.\ $10^{-4}$). HPR-QP solves the original Lasso~\eqref{eq:lasso-org}, others solve the CQP reformulation~\eqref{eq:lasso-convexqp}.  
‘T’ = time‐limit, ‘F’ = failure.}
\label{tab:uci-1e4}
\renewcommand{\arraystretch}{1.0}
\begin{tabularx}{\textwidth}{l *{5}{>{\centering\arraybackslash}X}}
\toprule
Instance             & HPR-QP & PDQP   & SCS & CuClarabel & Gurobi \\
\midrule
abalone7         & 1.4       & 78.3   & 1850   & 21.8       & 85.8   \\
bodyfat7         & 0.9       & 21.8   & T   & 1.6        & 27.5   \\
E2006.test       & 0.1       & 1.2    & 132.0   & 9.1        & 8.3    \\
E2006.train      & 0.3       & 1.7    & F   & 98.8       & 260.3  \\
housing7         & 3.9       & 44.7   & T   & 5.2        & 112.5  \\
log1p.E2006.test & 3.8       & 783.4  & T   & 170.0      & 120.4  \\
log1p.E2006.train& 11.7      & 1830.4 & T   & 310.0      & 551.2  \\
mpg7             & 0.1       & 6.5    & 411.0   & 0.2        & 1.1    \\
pyrim5           & 3.2       & 90.1   & T   & 3.2        & 33.7   \\
space\_ga9       & 0.3       & 13.6   & 332.0   & 5.5        & 31.0   \\
triazines4       & 57.0      & 1286.9 & T   & 25.3       & 310.1  \\
\midrule
SGM10 (Time)            & 4.3       & 74.1   & 1776.0   & 22.9       & 71.3   \\
\bottomrule
\end{tabularx}
\end{table}

\begin{table}[htp]
\centering
\caption{Numerical performance on randomly generated Lasso instances (Tol.\ $10^{-4}$).  
HPR-QP solves the original Lasso~\eqref{eq:lasso-org}, others solve the CQP reformulation~\eqref{eq:lasso-convexqp}.  
‘T’ = time‐limit, ‘M’ = out‐of‐memory.}
\label{tab:lasso-random-1e4}
\renewcommand{\arraystretch}{1.0}
\begin{tabularx}{\textwidth}{
  >{\centering\arraybackslash}p{1.5cm}  
  >{\centering\arraybackslash}p{1.5cm}  
  *{5}{>{\centering\arraybackslash}X}   
}
\toprule
$p$ & $q$ 
  & HPR-QP 
  & PDQP 
  & SCS 
  & CuClarabel 
  & Gurobi \\
\midrule
$10^4$         & $5\times10^5$  &  0.4  &  4.7  & T  & 2.1 & 21.9 \\
$10^4$         & $10^6$         &  0.1  &  5.5  & T  & 3.0 & 36.7 \\
$2\times10^5$  & $5\times10^6$  &  3.6  & 31.2  & T    & M &  T    \\
$2\times10^5$  & $10^7$         &  8.4  & 94.1  & T    & M &  T    \\
$4\times10^5$  & $5\times10^6$  &  5.5  & 49.5  & T    & M &  M    \\
$4\times10^5$  & $10^7$         & 16.2  &147.7  & T    & M &  M    \\
$6\times10^5$  & $5\times10^6$  &  8.1  & 62.6  & T    & M &  M    \\
$6\times10^5$  & $10^7$         & 27.9  &186.2  & T    & M &  M    \\
$8\times10^5$  & $5\times10^6$  & 10.6  & 78.0  & T    & M &  M    \\
$8\times10^5$  & $10^7$         & 31.8  &240.2  & T    & M &  M    \\
$10^6$         & $5\times10^6$  & 13.1  & 92.8  & T    & M &  M    \\
$10^6$         & $10^7$         & 46.1  &  M    & T   & M &  M    \\
\midrule
\multicolumn{2}{c}{{SGM10 (Time)}}            & 11.2 & 90.6  & 3600.0 &1395.2 &1684.4 \\
\bottomrule
\end{tabularx}
\end{table}

\begin{figure}[htp]
  \centering
  \begin{subfigure}[b]{0.48\textwidth}
    \centering
    \renewcommand{\arraystretch}{1.95}
    \begin{tabularx}{\linewidth}{l *{2}{>{\centering\arraybackslash}X}}
      \toprule
      Solver             & SGM10 (Time)  & Solved \\
      \midrule
      HPR-QP   &  0.4  & 36     \\
      PDQP               &  1.3  & 36     \\
      SCS                &  0.5  & 36     \\
      \midrule
      CuClarabel           &  4.2  & 36     \\
      Gurobi             & 22.8  & 36     \\
      \bottomrule
    \end{tabularx}
    \subcaption{SGM10 and number of problems solved}%
    \label{tab:real-QAP-1e4}
  \end{subfigure}%
  \hfill
  \begin{subfigure}[b]{0.48\textwidth}
    \centering
    \includegraphics[width=\linewidth]{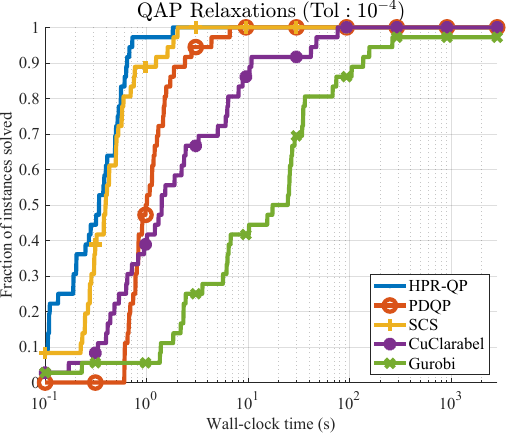}
    \subcaption{Absolute performance profiles}%
    \label{fig:real-QAP-1e4}
  \end{subfigure}

  \caption{Numerical performance of tested solvers on 36 QAP relaxations  (Tol.\ $10^{-4}$).}
  \label{fig:real-QAP-1e4-combined}
\end{figure}

\end{document}